\documentclass[a4paper]{article}
\pdfoutput=1
\makeatletter
\providecommand{\iftx@stxtwo}{\iffalse}
\providecommand{\iftx@ebgm}{\iffalse}
\providecommand{\iftx@ut}{\iffalse}
\providecommand{\iftx@nc}{\iftrue}
\providecommand{\iftx@fourier}{\iffalse}
\makeatother

\usepackage{amsthm}
\usepackage{tikz}
\usetikzlibrary{decorations.markings,decorations.pathreplacing,positioning,cd,fit,shapes.geometric}
\tikzset{radial vector/.style={blue, -latex}}
\tikzset{radial vector 2/.style={red, -latex}}
\tikzset{midarrow/.style={postaction=decorate,decoration={markings,mark={at position #1 with {\arrow{latex}}}}}}
\tikzset{invmidarrow/.style={postaction=decorate,decoration={markings,mark={between positions 0 and 1 step #1 with {\arrow{latex reversed}}}}}}
\usepackage{cancel}
\usepackage[top=2cm, bottom=2.5cm, outer=2cm, inner=2cm, marginparsep=0.7cm, marginparwidth=1.5cm]{geometry}
\usepackage[read,fonts=true,save=true,mdots=true,ptmx={all=upnewtx},lang=english,xtras=true,eseese=Terragnino,mfonts={lsym=ptm,cal=dutch,scr=true},bbold=true,textcomp=no]{mworksx}
\usepackage{newtxtextM,color}

\newcommand{\nipar}{\par\noindent}
\setint[2](2)
\makeatletter

\let\eatup\@gobble
\newcommand*\MapsTo{%
  \@ifstar{\xrightarrow[\raisebox{0.25 em}{\smash{\ensuremath{\sim}}}]{}}{\xrightarrow{\raisebox{-0.25 em}{\smash{\ensuremath{\sim}}}}}%
}
\newcommand{\fadecol}[1]{\textcolor{gray}{#1}}
\newif\ifmidv\midvtrue
\newcommand{\midv}{
    \newcommand{\fade}[1]{\fadecol{##1}}
    \newcommand{\fadees}{\fade}
    \newcommand{\fadecap}[3]{\fadees{##1}{##2}{\MakeUppercase##3}}
    \newcommand{\lac}[1]{\@ifstar{<>##1}{<##1>}}
    \newcommand{\slac}[1]{<##1> }
    \newcommand{\xlac}[1]{\textcolor{red}{<}##1\textcolor{red}{>}\@ifstar{\ }{}}
    \newcommand{\xsqb}{\sqb}
    \midvtrue
}
\newcommand{\finv}{
    \newcommand{\fade}[1]{}
    \newcommand{\fadees}{\fade}
    \newcommand{\fadecap}[3]{\fadees{##1}{##2}{\MakeUppercase##3}}
    \newcommand{\lac}{\@ifstar{\expandafter\@gobble\@gobble}{}}
    \newcommand{\slac}{\lac}
    \newcommand{\xlac}{\@ifstar{\lac}{\lac}}
    \newcommand{\xsqb}{}
    \finvtrue
}
\makeatother
\upcaps

\newcommand{\opn}{\operatorname}

\setdiff[0](.25)

\newcommand{\vq}{\vartheta}

\newcommand{\lsim}{\lesssim}
\midv
\renewcommand{\thesection}{\arabic{section}}

\makeatletter
\newcommand{\frl}{\@ifstar{\@frl}{\frl@}}
\newcommand{\@frl}[1]{(-\Dg)^{\xfr{#1}{2}}}
\newcommand{\frl@}[1]{(-\Dg)^{#1}}
\makeatother
\definethm{probl}{Problem}[s]
\definethm{psprobl}{Pseudo-problem}[s]
\definethm{typo}{Probable typo/misprint}[s]
\definethm{genoss}{Navier-Stokes vs. Euler Remark}[s]

\makeatletter
\newcommand{\xtag}{\@ifstar{\@tag}{\tag@}}
\newcommand{\@tag}{\stepcounter{equation}\tag{\theequation}}
\newcommand{\tag@}[1]{\@tag}
\newcommand{\stag}[1]{\@tag\label{eq:#1}}
\makeatother
\numberwithin{equation}{section}
\makeatletter
\newcommand{\xeqref}[2]{
\expandafter\ifx\csname r@#2\endcsname\relax
  \protect\G@refundefinedtrue
  \nfss@text{\reset@font\bfseries (??) R#1}%
  \@latex@warning{%
      Reference `#2' on page \thepage \space undefined%
  }%
\else
  \expandafter\expandafter\expandafter\Hy@setref@link\csname r@#2\endcsname\@empty\@empty\@nil{{\textup{\tagform@{\ref*{#2}}} R#1}\expandafter\@gobble\@firstoffive}
  \fi}
\newcommand{\EXP}{\qg}
\newcommand{\xfrl}{\@ifstar{\frl*\EXP}{\frl\EXP}}
\makeatother

\let\dacute\H
\def\H{\mathbb H}

\def\L{\mathbb L}

\newcommand{\pon}{\opn}
\providecommand{\Cinf}{\s{C}^\8}

\title{Density of weak solutions of the fractional Navier-Stokes equations in the smooth divergence-free vector fields}
\author{Michele Gorini}
\date{}

\begin{document}
\maketitle
\begin{abstract}
In this paper, we consider the fractional Navier-Stokes equations. We extend a previous non-uniqueness result due to Cheskidov and Luo, found in \cite{CL}, from Navier-Stokes to the fractional case, and from $L^1$-in-time, $W^{1,q}$-in-space solutions for every $q>1$ to $L^s$-in-time, $W^{1,q}$-in-space solutions for appropriate ranges of $s,q$.
\end{abstract}
\sect{Introduction}
\ssect{Background}
\ssectcnt*{teor}\ssectcnt*{propo}\ssectcnt*{defi}\ssectcnt*{oss}\ssectcnt*{cor}\ssectcnt*{lemma}
We consider the fractional Navier-Stokes equations on the $d$-dimensional torus $\T^d=[-\pg,\pg]^d$:
\[\begin{sistema}
\pd_tv+\opn{div}(v\otimes v)+\grad p+\xfrl v=0 \\
\opn{div}v=0
\end{sistema}, \stag{FNS}\]
where $v=(v_1,v_2,v_3)^T(t,x)\in\R^d$ is the velocity and $p(t,x)\in\R$ is the pressure of the fluid, and the fractional laplacian, for $\EXP\geq0$, is defined vie the Fourier transform:
\[\xfrl u(t,k)=|k|^{2\EXP}\hat u(t,k).\]
For $\EXP=1$, these are the Navier-Stokes equations, a fundamental mathematical model describing the motion of an incompressible viscous Newtonian fluid. The fractional case $\EXP\neq1$, for suitable exponents ($\EXP\in[\xfr12,1]$), can also be used to model fluid flow, as described in \cite{FNSFis}.
\nipar The corresponding model for inviscid fluids is given by the Euler equations:
\[\begin{sistema}
\pd_tv+\opn{div}(v\otimes v)+\grad p=0 \\
\opn{div}v=0
\end{sistema}. \stag{Eul}\]
We consider weak solutions defined as space-time distributional solutions.
\xbegin{defi}[Weak solutions]
Let $\s{D}_T$ be the space of divergence-free test functions $\fg\in\Cinf(\R\x\T^d)$ such that $\fg=0$ if $t\geq T$. Let $u_0\in L^2(\T^d)$ be weakly divergence-free. A vector field $u\in L^2_tL^2_x([0,T]\x\T^d)$ is a \emph{weak solution of \eqref{eq:FNS} with initial data $u_0$} if the following hold:
\ben[label={\upshape\arabic*)}]
\item For a.e. $t\in[0,T]$, $u(t,\per)$ is weakly divergence-free;
\item For any $\fg\in\s{D}_T$,
\[\xints{\T^d}{}\2u_0(x)\per\fg(0,x)\diff x=-\xints0T\xints{\T^d}{}\2u\per(\pd_t\fg-\xfrl\fg+u\per\grad\fg)\diff x\diff t.\]
\een
\emph{Weak solutions for the Euler equations \eqref{eq:Eul}} are defined similarly, by removing the term $u\per\xfrl\fg$ from the integral in the formula above.
\xend{defi}
A more physical class of solutions, which Leray introduced in the case $\EXP=1$ and proved to exist in $\R^3$ in \cite{L}, and Hopf proved to exist in general domains in $d\geq2$ for $\EXP=1$, is that of so-called ``admissible'' or ``Leray-Hopf'' weak solutions, which satisfy an energy inequality.
\xbegin{defi}[Admissible solutions]
A weak solution of \eqref{eq:FNS} is called an \emph{admissible weak solution} if $u\in\s{C}_w([0,T];L^2)\cap L^2_tH^\EXP_x([0,T]\x\T^d)$ and
\[\xfr12\norm{u(t)}_{L^2(\T^d)}^2+\xints0t\norm{\xfrl*u(s)}_{L^2(\T^d)}^2\diff s\leq\xfr12\norm{u(0)}_{L^2(\T^d)}^2\]
for all $t\in[0,T]$.
\nipar A similar definition is given for \eqref{eq:Eul}, once again suppressing the $\xfrl$ term from the energy inequality, and not requiring the $L^2_tH^\EXP_x$ regularity.
\end{defi}
\xbegin{oss}[Euler and Leray]
The term ``Leray-Hopf weak solution'' is not used for the Euler equations because there is no existence result for admissible solutions of the Euler equations with generic $L^2$ initial data.
\xend{oss}
This energy inequality, also called admissibility condition, is a relaxation of the natural conservation law of \eqref{eq:FNS} (or \eqref{eq:Eul}). Solutions which satisfy it enjoy much better properties than general weak solutions. At least in the case $\EXP=1$, there is a vast literature on the topic, giving results such as weak-strong uniqueness \cite{P,S,La} and partial regularity \cite{La,S2}. For more such properties, cfr. \cite{CL,LQZZh} and the references therein. Even in the fractional case, it is relatively simple to prove weak-strong uniqueness for $\s{C}^1$ solutions.
\nipar Such properties, though nice from a regularity point of view, make it very hard to construct non-unique Leray-Hopf weak solutions in $d\geq3$(\fn{2D solutions in the case $\EXP=1$ are known to be unique and smooth.}). Indeed, the problem of uniqueness (or non-uniqueness) of Leray-Hopf weak solutions remains a challenging open one. So far, there are only numerical evidence \cite{NENU} and partial results.
\nipar One direction of these partial results is to work with forced Navier-Stokes, namely the following system:
\[\begin{sistema}
\pd_tv+\opn{div}(v\otimes v)+\grad p+\xfrl v=f\in L^1_tL^2_x \\
\opn{div}v=0
\end{sistema}, \stag{ForcNS}\]
with admissibility condition
\[\xfr12\norm{u(t)}_{L^2(\T^d)}^2+\xints0t\norm{\xfrl*u(s)}_{L^2(\T^d)}^2\diff s\leq\xfr12\norm{u(0)}_{L^2(\T^d)}^2+\xints0t\xints{\T^3}{}\2f\per v\diff x\diff s\]
A non-uniqueness result for Leray-Hopf weak solutions of \eqref{eq:ForcNS}, for a suitable forcing term $f$, was proved in \cite{ABC}.
\nipar The other direction is to consider general solutions which may or may not have sufficient regularity to be tested for admissibility. The first result in this direction is \cite[Theorem 1.2]{BV}.
\xbegin{teor}[Buckmaster-Vicol][thm:teor:BV]
There exists $\bg>0$ such that, for any non-negative smooth function $e(t):[0,T]\to\R_{\geq0}$, there exists $v\in\s{C}^0_tH^\bg_x([0,T]\x\T^3)$ a weak solution of the Navier-Stokes equations such that, for all $t\in[0,T]$,
\[e(t)=\xints{\T^3}{}\2\abs{v(t,x)}^2\diff x.\]
Moreover, the associated vorticity $\grad\x v$ lies in $\s{C}^0_tL^1_x([0,T]\x\T^3)$.
\xend{teor}
The result of \cite{BV} immediately yields non-uniqueness for zero initial data, but nontrivial solutions with zero initial data cannot be admissible, since the kinetic energy itself has to increase.
\nipar A similar result was obtained for the Euler equations \eqref{eq:Eul} in \cite{BDLSzV} for $\s{C}^{\xfr13-}$ solutions.
\xbegin{teor}[Buckmaster-De Lellis-Székelyhidi-Vicol][thm:teor:BDLSzV]
For every $\bg<{}^1\!/_3$ and every positive smooth $E:[0,T]\to\R$, there exists a solution $(v,p)\in\s{C}^\bg([0,T]\x\T^3)$ of the Euler equations such that
\[\xfr12\xints{\T^3}{}\2\abs{v(t,x)}^2\diff x=E(t).\]
\xend{teor}
This was in the context of the Onsager conjecture (formulated in \cite{O} by Lars Onsager), which states that admissible solutions of the Euler equations (i.e. which satisfy the energy inequality for those equations, which reads $\|u(t,\per)\|_{L^2}\leq\|u(0,\per)\|_{L^2}$) are unique in $\s{C}^\bg$ for $\bg>\xfr13$, whereas non-uniqueness holds for $\bg<\xfr13$.
\nipar The presence of the Laplacian term in the fractional Navier-Stokes equations is an obstacle in proving non-uniqueness of Leray-Hopf weak solutions, as it proves hard to control. However, if the exponent is small enough ($\EXP<\xfr13$), non-uniqueness can still be proved in a very strong fashion (cfr. \cite{CDLDR,DR,G}). For $\EXP\geq\xfr13$, we are led back to general solutions.
\nipar Another non-uniqueness result was proved in \cite{BCV}, both for $\EXP=1$ and for some fractional values of $\EXP$. The authors of the paper first proved the following gluing theorem, which is \cite[Theorem 1.5]{BCV}, and from there deduced the non-uniqueness of $\s{C}^0_tH^\bg_x$ weak solutions of \eqref{eq:FNS}, for any $\ag\in[1,\xfr54)$, any initial datum $v_0\in\dot H^3$, and any sufficiently small $\bg$.
\xbegin{teor}[Buckmaster-Colombo-Vicol][thm:teor:BCV]
For $\EXP\in[1,\xfr54)$, there exists a $\bg=\bg(\EXP)>0$(\fn{The maximal $\bg_{max}(1)$ for which this holds in the case $\EXP=1$ can be quantified as $\bg_{max}\dsim10^{-3}$.}) such that the following holds. For $T>0$, let $v\stp 1,v\stp 2\in\s{C}^0([0,T];\dot H^3(\T^3)\!)$ be two strong solutions of the Navier-Stokes equations \eqref{eq:FNS} on $[0,T]$, with data $v\stp 1(0,x)$ and $v\stp 2(0,x)$ of zero mean. There exists a weak solution $\lbar v$ of the Cauchy problem to \eqref{eq:FNS} on $[0,T]$ with initial datum $\lbar v|_{t=0}=u\stp 1|_{t=0}$, which has the additional regularity $\lbar v\in\s{C}^0([0,T];H^\bg(\T^3)\cap W^{1,1+\bg}(\T^3)\!)$, and such that $\lbar v\equiv v\stp 1$ on $[0,{}^T\!/{\!}_3]$ and $\lbar v\equiv v\stp 2$ on $[{}^2\!/{\!}_3T,T]$. Moreover, for every such $\lbar v$, there exists a zero Lebesgue measure set of times $\Sg_T\sbs[0,T]$ with Hausdorff (in fact box-counting) dimension less than $1-\bg$ such that $\lbar v\in\Cinf(\!(\!(0,T]\ssm\Sg_T)\x\T^3)$. In particular, $\lbar v$ is smooth almost everywhere.
\xend{teor}
Returning to $\EXP=1$, in \cite{CL}, Cheskidov and Luo prove non-uniqueness for $L^1_tW^{1,q}_x$ solutions for any $q$. They do this via a density theorem, \cite[Theorem 1.7]{CL}, which I will frame as a ``meta-theorem'' plus a choice of parameter ranges.
\xbegin{teor}[Meta-theorem]
Let $d\geq2$ be the dimension, $\EXP>0$, $1\leq p<2,q,s<\8$, and $\gg,\eg>0$. For any smooth, divergence-free vector field $v\in\Cinf([0,T]\x\T^d)$ with zero spatial mean for each $t\in[0,T]$, henceforth known as ``starting field'', there exists a weak solution $u$ of \eqref{eq:FNS} and a set:
\[\s{I}=\bigcup_{i=1}^\8(a_i,b_i)\sbs[0,T],\]
such that the following holds.
\ben[label={\upshape(\arabic*)}]
\item The solution u satisfies $u\in L^p_tL^\8_x([0,T]\x\T^d)\cap L^s_tW^{\gg,q}_x([0,T]\x\T^d)$;
\item $u$ is a smooth solution on $(a_i,b_i)$ for every $i$. Namely, $u|_{\s{I}\x\T^d}\in\Cinf(\s{I}\x\T^d)$. In addition, $u$ agrees with the smooth solution emerging from the initial data $v(0)$ near $t=0$; in fact, if there is an interval $[a,b]\sbse[0,T]$ where $v$ is an exact solution, then $u=v$ on $[a,b]$;
\item The Hausdorff dimension of the residue set $\s{S}=[0,T]\ssm\s{I}$ satisfies $d_{\s{H}}\s{S})\leq\eg$;
\item The solution $u$ and the starting field $v$ are $\eg$-close in $L^p_tL^\8_x\cap L^1_tW^{1,q}_x$.
\een
\xend{teor}
Taking any smooth datum $v_0$, one can consider the smooth solution $\tilde v$ with that initial datum, and glue it to a different smooth field $u$ by use of a cutoff function:
\[\tilde u\coloneq\xg\tilde v+(1-\xg)u,\]
where $\xg$ is 1 near $t=0$ and $\xg=1$ near $t=T$. By applying the meta-theorem to $\tilde u$, one quickly concludes both non-uniqueness of solutions and a gluing theorem à la \cite{BCV} for the ranges where the meta-theorem can be proved to hold.
\nipar In this framing, \cite[Theorem 1.7]{CL} can be stated as follows.
\xbegin{teor}[Cheskidov-Luo][thm:teor:CL]
The meta-theorem holds for $\EXP=1,s=1,\gg=0$, and any $q$.
\xend{teor}
Compared to \cite{BCV}, this suggests there may be a tradeoff between space regularity and time regularity: giving up some time regularity ($\s{C}^0$ in \cite{BCV} vs. $L^1$ in \cite{CL}), we have gained space regularity ($H^\bg$ with small $\bg$ in \cite{BCV} vs. $W^{1,q}$ for any $q$ in \cite{CL}).
\nipar This result was then extended to $L^s_tW^{\gg,q}_x$ and $\EXP\neq1$ in \cite[Theorem 1.2]{LQZZh}, which can be restated as follows.
\xbegin{teor}[Li-Qu-Zong-Zhang][thm:teor:LQZZh]
The meta-theorem holds whenever one of the following holds:
\begin{align*}
(\EXP,\gg,s,q)\in{}&\sqpa{\xfr54,2}\x[0,3)\x[0,\8]\x[1,\8]\quad0\leq\gg<\xfr{4\EXP-5}{s}+\xfr3q+1-2\EXP \\
(\EXP,\gg,s,q)\in{}&[1,2)\x[0,3)\x[1,\8]\x[1,\8]\quad0\leq\gg<\xfr{2\EXP}{s}+\xfr{2\EXP-2}{q}+1-2\EXP.
\end{align*}
The solutions can also be chosen in $H^{\bg'}_{t,x}$ for sufficiently small $\bg'<1$, and with supports close to that of the starting field.
\end{teor}
\clearpage

\ssect{Main theorem}
This paper is devoted to proving the following theorem and corollary.
\xbegin{teor}[][thm:teor:Relaz]
Let $d\geq2$ be the dimension, $\eg>0$, and $p,q,s,\bg',\dg',\zg',\gg,\EXP$ satisfy the following system:
\begin{align*}
\bg',\dg',\zg'>{}&0 \\
p,s\geq{}&1 \\
q>{}&1 \\
\dg'<{}&2\bg' \\
\xfr{p}{2-p}(d-1)<{}&\zg' \\
s\pa{\dg'\gg+\xfr{\zg'}{2}+\gg+\xfr{d-1}{2}-\xfr{d-1}{q}}<{}&\zg' \\
2\bg'+\zg'<{}&2\dg'+d+1 \\
2\EXP\leq{}&\gg+1 \\
\gg\leq{}&d \\
q\leq{}&\xfr{2d}{2\gg-d}\one_{2\gg>d}+\8\one_{2\gg\leq d}. \\
\end{align*}
Then, the meta-theorem holds for $p,q,s,\EXP$.
\xend{teor}
By studying the system, one can arrive at the following definitions and corollary.
\begin{align*}
Q(d,\gg)\coloneq{}&\xfr{d-1}{\gg-1} \\
q_0(\gg,d)\coloneq{}&\xfr{2d-2}{d-1+2\gg}<q_G \\
S_1(p,d,q,\gg)\coloneq{}&\xfr{p(d-1)}{d-1+(2-p)\pa{\gg-\xfr{d-1}{q}}} \\
S_2(d,q,\gg)\coloneq{}&\xfr{d+1}{d+\gg-\xfr{d-1}{q}} \\
S_3(\gg,d,q)\coloneq{}&\xfr{2\gg+d-1}{2\gg+d-1-\xfr{d-1}{q}}
\end{align*}
\xbegin{cor}[][thm:cor:Ranges]
Assume $p<1+\xfr1d$ and $2\EXP\leq\gg+1$. The meta-theorem holds:
\bi
\item If $\gg<1,q>q_0,s<S_2$;
\item If $\gg<1,p<1+\xfr{\gg}{\gg+d-1},1<q\leq q_0,s<S_3$;
\item If $\gg<1,p\geq1+\xfr{\gg}{\gg+d-1},q\leq q_0,s<S_1\leq S_1(1,d,q,\gg)=\xfr{d-1}{d-1+\gg-\xfr{d-1}{q}}$;
\item If $\gg\in[1,d),q\in(q_0,Q),s<S_2$;
\item If $\gg\in[1,\xfr{d-1}{2}),q\in(1,q_0],s<S_1\leq S_1(1,d,q,\gg)=\xfr{d-1}{d-1+\gg-\xfr{d-1}{q}}$.
\ei
\xend{cor}
The paper is organized as follows. In section \ref{Outl}, we give an outline of the proof, stating the main iteration proposition, which lies at the heart of the proof, at the end of the section. In section \ref{Gl}, we carry out the first substep of the convex integration step, which is a gluing argument. In section \ref{Conc}, we show how the second substep is done, namely the perturbation step. All three of these sections are very similar to the corresponding sections of \cite{CL}, so the various lemmas and propositions are mostly left without proof, as the proofs can be found in \cite{CL}. In section \ref{Est}, we show that, given a set of relations between the parameters of the convex integration step, we can obtain the satisfactory estimates that lead to the proof of the meta-theorem. This section is also quite similar to the corresponding one in \cite{CL}, except that, instead of choosing the parameters first and immediately deducing the needed relations, we postpone the choice of the parameters to section \ref{Calcs}. In that last section, we gather the relations, adding an extra couple which arise in the gluing argument of section \ref{Gl}, and study the resulting system, completing the proofs of \kcref{thm:teor:Relaz} and \kcref{thm:cor:Ranges}.

\sect{Outline of the proof}\label{Outl}
The proof of \kcref{thm:teor:Relaz} consists of an iterative scheme achieved by  a repeated application of this paper's main proposition, \kcref{thm:propo:MPS}, to obtain a sequence of solutions to the fractional Navier-Stokes-Reynolds system below:
\[\begin{sistema}
\pd_tv+\opn{div}(v\otimes v)+\grad p+\xfrl v=-\opn{div}R \\
\opn{div}v=0
\end{sistema}.
\stag{FNSR}.\]
We are omitting the pressure when referring to solutions of \eqref{eq:FNSR} because it can be uniquely determined by the following elliptic equation, provided it is average-free:
\[\D p=\opn{div}\opn{div}R-\opn{div}\opn{div}(v\otimes v).\]
The proof mainly consists of three goals:
\ben[label={\upshape(\alph*)}]
\item The convergence of $u_n\to u$ in $L^2_tL^2_x$ and of $R_n\to0$ in $L^1_tL^1_x$ so that $u$ is a weak solutions of \eqref{eq:FNS};
\item The convergence of $u_n\to u$ in $L^p_tL^q_x\cap L^s_tW^{\gg,q}_x$;
\item The small dimension of the set of the singular times of $u$.
\een
To this end, as done in \cite{CL}, we employ a two-step approach:
\bi
\item Step 1: $(u_n,R_n)$ is transformed into $(\lbar u_n,\lbar R_n)$ by concentrating the stress tensor;
\item Step 2: space-time convex integration turning $(\lbar u_n,\lbar R_n)$ into $(u{n+1},R_{n+1})$.
\ei
The first step is mainly to achieve a small singular set in time, and the second step ensures the convergences.

\ssect{Step 1: Concentrating the stress error}
\ssectcnt*{teor}\ssectcnt*{propo}\ssectcnt*{defi}\ssectcnt*{oss}\ssectcnt*{cor}\ssectcnt*{lemma}
Given $(u_{n-1},R_{n-1})$, we divide the time interval $[0,T]$ into smaller sub-intervals $I_i$ of length $\tg^\eg$, where $\tg>0$ will be chosen to vary depending on $(u_{n-1},R_{n-1})$. The total number of sub-intervals is thus of order $\tg^{-\eg}$.
\nipar On each $I_i$, we solve a difference fractional Navier-Stokes system centered at $(u_{n-1},R_{n-1})$ to obtain a corrector $v_i$ on $I_i$. More precisely, $v_i:I_i\x\T^d\to\R^d$ solves
\[\begin{sistema}
\pd_tv_i-\Dg v_i+\opn{div}(v_i\otimes v_i)+\opn{div}(v_i\otimes u)+\opn{div}(u\otimes v_i)+\grad q_i=\opn{div}R \\
\opn{div}v_i=0 \\
v_i(t_i)=0
\end{sistema},\]
so that $u_{n-1}+v_i$ is an exact solution of the fractional Navier-Stokes equations \eqref{eq:FNS} on $I_i$.
\nipar To concentrate the error and obtain a solution of \eqref{eq:FNSR} on $[0,T]$, we apply a sharp cutoff $\xg_i$ to the corrector $v_i$ and obtain the glued solution $\lbar u_{n-1}$ defined by:
\[\lbar u_{n-1}\coloneq u_{n-1}+\sum_i\xg_iv_i.\]
Specifically, each $\xg_i$ equals 1 on a majority of $I_i$, but $\xg_i\equiv0$ near endpoint regionss of $I_i$ of length of order $\tg$. Since $\eg\ll1$, the cutoff $\xg_i$ is very sharp compared to the length of the sub-interval $I_i$.
\nipar On one hand, due to the sharp cutoff $\xg_i$, the stress error $\lbar R_{n-1}$ associated with $\lbar u_{n-1}$ will only be supported on endpoint regions of $I_i$ of length of order $\tg$. In other words, the temporal support of $\lbar u_n$ can be covered by $\sim\tg^{-\eg}$ many intervals of size $\sim\tg$, from which one already sees the singular set of the final solution will have a small dimension.
\nipar On the other hand, the corrector $v_i$ is very small, say in $L^\8_tH^d_x$, since it starts with initial data 0 and we can choose time scale $\tg^\eg=|I_i|$ to be sufficiently small. More importantly, the new stress error $\lbar R_{n-1}$ associated with $\lbar u_{n-1}$ satisfies the estimate:
\[\norm{\lbar R_{n-1}}_{L^1_tL^r_x}\lsim\norm{R_{n-1}}_{L^1_tL^r_x}\qquad1<r<\8,\]
with an implicit constant independent of the time scale $\tg>0$. In other words, concentrating the stress error $R_{n-1}$ to $\lbar R_{n-1}$ cost a loss of a constant multiple when measuring in $L^1$ norm in time.

\ssect{Step 2: Space-time convex integration}
The next step is to use a convex integration technique to reduce the size of $\lbar R_{n-1}$ by adding a perturbation $w_n$ to $\lbar u_{n-1}$ to obtain a new solution $(u_n,R_n)$ of \eqref{eq:FNSR}. The perturbation $w_n$ and the new stress $R_n$ satisfy the equation
\[\opn{div}R_n=\opn{div}\lbar R_{n-1}+\opn{div}(w_n\otimes w_n)+\pd_tw_n-\Dg w_n+\opn{div}(\lbar u_{n-1}\otimes w_n)+\opn{div}(w_n\otimes\lbar u_{n-1})+\grad P_n,\]
for a suitable pressure $P_n$. Heuristically, we wish to balance the old Reynolds stress with the quadratic term $w_n\otimes w_n$, that is:
\[\opn{div}(\lbar R_{n-1}+w_n\otimes w_n)=HSF+HTF+LO. \stag{Bal}\]
Here, HSF are terms of high spatial frequency, HTF have high temporal frequency and will be balanced further by a part of $\pd_tw_n$, and LO are lower-order terms. This is similar to \cite{BV,BCV}, but with the fundamental difference that this additional ``convex integration in time'' requires no additional oscillation/concentration constraint and is basically free, which is crucial in obtaining the regularity ranges $L^p_tL^q_t\cap L^s_tW^{\gg,q}_x$ of our main theorem.
\nipar In order to obtain such balances, we need two ingredients:
\ben[label=(\arabic*)]
\item Suitable stationary flows as the spatial building blocks of our perturbation; these must be able to achieve some level of spatial concentration;
\item Intermittent temporal functions to oscillate the spatial building blocks in time.
\een
Once the first ingredient is found, the second one is relatively straightforward to implement. On the technical side, the first ingredient will be the Mikado flows introduced in \cite{DSz}. These are periodic pipe flows that can be arranged to be supported on periodic cylinders with small radius. In other words, Mikado flows can achieve $(d-1)$-dimensional concentration on $\T^d$. It is worth noting that, in the framework of \cite{BV}, stationary Mikado flows are not sufficiently intermittent to be used for the Navier-Stokes equations (\eqref{eq:FNS} with $\EXP=1$) in dimension $d\leq3$.
\nipar The desired balance \eqref{eq:Bal} imposes a relation between $w_n$ and $\lbar R_{n-1}$:
\[\norm{w_n}_{L^2_{t,x}}\sim\norm{\lbar R_{n-1}}_{L^1_{t,x}}. \stag{wnRn-1}\]
This relation will imply the $L^2_{t,x}$ convergence of the approximate solutions $w_n$ as long as one can successfully reduce the size of the stress error:
\[\norm{R_n}_{L^1_{t,x}}\ll\norm{\lbar R_{n-1}}_{L^1_{t,x}}.\]
In particular, special attention will be paid to estimating the temporal derivative component of the new Reynolds stress, defined by
\[\opn{div}R_{tem}=\pd_tw_n, \stag{Rtempder}\]
and achieving the regularity of the perturbation:
\[\norm{w_n}_{L^p_tL^q_x}+\norm{w_n}_{L^s_tW^{\gg,q}_x}\ll1. \stag{PertReg}\]
These two constraints \eqref{eq:Rtempder} and \eqref{eq:PertReg} require a very delicate parameter choice when designing the perturbation. On one hand, \eqref{eq:Rtempder} implies the temporal frequency cannot bee too large, relative to the spatial frequency, otherwise the time derivative will dominate. On the other hand, \eqref{eq:PertReg} requires a large temporal frequency, so as to use the temporal concentration to offset the loss caused by going from $L^2$ to $L^q$ or $W^{\gg,q}$ in space in relation to \eqref{eq:wnRn-1}. This is achieved by, roughly speaking, a tradeoff: we trade $L^2$ in time, reducing it to $L^p$ or $L^s$, in order to obtain $L^q$ or $W^{\gg,q}$ in space.

\ssect{Oscillation and concentration}
We do this computation in general dimension $d\geq2$ and $D\in[0,d]$ denotes the spatial intermittency.
\nipar We start with a velocity perturbatoin in $L^2_{t,x}$ with a certain decay given by the previous stress error:
\[\norm{w_n}_{L^2_{t,x}}\lsim1.\]
Denote the spatial frequency by $\lg$ and the temporal frequency by $\kg$, namely:
\[\norm{\pd_t^m\grad^nw_n}_{L^2_{t,x}}\lsim\kg^m\lg^n.\]
The intermittency parameter $D$ in space dictates the level of concentration of $w_n$ and the scaling law:
\[\norm{w_n(t)}_{L^q}\lsim\norm{w_n(t)}_{L^2}\lg^{(d-D)\pa{\xfr12-\xfr1q}}.\]
As for the temporal scaling, we assume for simplicity that $w_n$ is fully concentrated in time:
\[\norm{w_n}_{L^p_tL^q_x}\lsim\norm{w_n}_{L^2_tL^q_x}\lg^{\xfr12-\xfr1p}\sim\kg^{\xfr12-\xfr1p}\lg^{(d-D)\pa{\xfr12-\xfr1q}}.\]
With such scaling laws, we effectively assume negligible temporal oscillation and the goal then reduces to finding a choice of $D$ in terms of the parameters $d,p,q$. In other words, we need to find a balance between spatial oscillation and spatial concentration.
\nipar By the scaling relations just above, the stress error contributed by the time derivative \eqref{eq:Rtempder} satisfies:
\[\norm{\opn{div}^{-1}\pd_tw_n}_{L^1_{t,x}}\lsim\kg^{\xfr12}\lg^{-1}\lg^{\xfr{d-D}{2}}, \stag{RtempderEst}\]
where we assume $\opn{div}^{-1}$ gains one full derivative in space. The regularity condition \eqref{eq:PertReg} then becomes:
\begin{align*}
\norm{w_n}_{L^p_tL^q_x}\sim\kg^{\xfr12-\xfr1p}\lg^{\xfr{d-D}{2}}\ll1 \\
\norm{w_n}_{L^s_tW^{\gg,q}_x}\sim\kg^{\xfr12-\xfr1s}\lg^{1+\gg+\xfr{d-D}{2}-\xfr{d-D}{q}}. \stag{RegWEst}
\end{align*}
Conditions \eqref{eq:RtempderEst} and \eqref{eq:RegWEst} imply that:
\[\kg^{\xfr1s}\lg^{1+\gg+\xfr{d-D}{2}-\xfr{d-D}{q}}\ll\kg^{\xfr12}\ll\lg^{1+\xfr{d-D}{2}}.\]
Solutions to this exist for the ranges of parameters described in Section \ref{Calcs}.

\ssect{The main iteration proposition}
We are ready to introduce the main iteration proposition of the paper that materializes the above discussion. To simplify presentation, let us introduce the notion of well-prepared solutions to \eqref{eq:FNSR}, which encodes the small Hausdorff dimension of the set of singular times. Throughout the paper, we take $T=1$ and assume $0<\eg<1$ without loss of generality.
\xbegin{defi}[Well-prepared solution][thm:defi:WPS]
We say a smooth solution $(u,R)$ of \eqref{eq:FNSR} on $[0,1]$ is \emph{well-prepared} if there exist a set $I$ and a length scale $\tg>0$ such that $I$ is a union of at most $\tg^{-\eg}$-many closed intervals of length $5\tg$, and:
\[R(t,x)=0\qquad\VA t:\opn{dist}(t,I^c)\leq\tg.\]
\xend{defi}
With this definition, to ensure the solution $u$ has intervals of regularity with a small residue set of Hausdorff dimension $\lsim\eg$, it suffices to construct approximate solutions $(u_n,R_n)$ that are well-prepared for some $I_n,\tg_n$ such that:
\[I_n\sbs I_{n-1}\qquad\tg_n\to0.\]
The main proposition of this paper states as follows.
\xbegin{propo}[Main iteration][thm:propo:MI]
For any $\eg>0$, there exists a universal constant $M=M(\eg)>0$ such that for any $p,q,s,\gg,\EXP$ in the ranges of \kcref{thm:teor:Relaz} and $d\geq2$ there exists $r=r(p,q,s,d,\gg,\EXP)>1$ such that the following holds.
\nipar Let $\dg>0$ and $(u,R)$ be a well-prepared smooth solution of \eqref{eq:FNSR} for some set $\tilde I$ and a length scale $\tilde\tg>0$. Then there exists another well-prepared smooth solution $(u_1,R_1)$ of \eqref{eq:FNSR} for some set $I\sbs\tilde I$ with $0,1\nin I$, and some time scale $\tg<{}^{\tilde\tg}\!/_2$, such that:
\[\norm{R_1}_{L^1_tL^r_x}\leq\dg.\]
Moreover, the velocity perturbation $w\coloneq u_1-u$ satisfies:
\begin{align}
\opn{supp}w\sbs{}&I\x\T^d \\
\norm{w}_{L^2_tL^2_x}\leq{}&M\norm{\lbar R}_{L^1_tL^1_x} \label{MI:L2L2} \\
\norm{w}_{L^p_tL^q_x}+\norm w_{L^s_tW^{\gg,q}_x}\leq{}&\dg. \label{MI:LpLq+LsW1q}
\end{align}
\xend{propo}
A couple of comments to close this section:
\bi
\item The parameter $r>1$ is used to ensure the $L^r$ boundedness of the Calderón-Zygmund singular integral and is very close to 1;
\item Due to the local well-preparedness, on large portions of the time axis the solutions are exact solutions of the fractional Navier-Stokes equations \eqref{eq:FNS}, and we do not touch them in the future.
\ei
The main theorem \kcref{thm:teor:Relaz} is deduced from this proposition with precisely the same arguments of \cite[Section 2.6, proof of Theorem 1.7]{CL}.

\sect{Concentrating the stress error} \label{Gl}
The goal of this section is to prove \kcref{thm:propo:StressConc} below. The idea is that, given a solution $(u,R)$ of \eqref{eq:FNSR}, we can add a small correction term to it so that all of the stress error $R$ concentrates on a set $I$, the union of small intervals of length $\tg$, and thus obtain a new solution $(\lbar u,\lbar R)$. The key is that the procedure $(u,R)\leadsto(\lbar u,\lbar R)$ leaves the size of the stress $R$ invariant in $L^1_t$, up to a cost of a constant multiple:
\[\norm{\lbar R}_{L^1_tL^r_x}\leq C\norm{R}_{L^1_tL^r_x}\qquad r\in(1,\8),\]
where $C=C(r,\eg)$ is a universal constant that only depends on the exponent $r$ and the well-preparedness parameter $\eg>0$.
\xbegin{propo}[Error concentration][thm:propo:StressConc]
Let $0<\eg<1$ and $(u,R)$ be a well-prepared smooth solution of \eqref{eq:FNSR} for some set $\tilde I$ and length scale $\tilde\tg>0$. For any $1<r<\8$, there exists a universal constant $C=C(r,\eg)>0$ such that the following holds.
\nipar For any $\dg>0$, there exists another well-prepared smooth solution $(\lbar u,\lbar R)$ of \eqref{eq:FNSR}, for some set $0,1\nin I\sbs\tilde I$ and length scale $\tg<{}^{\tilde\tg}\!/_2$, satisfying the following:
\ben[label={\upshape(\arabic*)}]
\item The new stress error $\lbar R$ satisfies:
\begin{align*}
\lbar R(t,x)={}&0\qquad\opn{dist}(t,I^C)\leq\xfr32\tg \\
\norm{\lbar R}_{L^1_tL^r_x([0,1]\x\T^d)}\leq{}&C\norm{R}_{L^1_tL^r_x([0,1]\x\T^d)};
\end{align*}
\item The velocity perturbation $\lbar w\coloneq\lbar u-u$ satisfies:
\begin{align*}
\opn{supp}\lbar w\sbs{}&\tilde I\x\T^d \\
\norm{\lbar w}_{L^\8_tH^d_x([0,1]\x\T^d)}\leq\dg.
\end{align*}
\een
\xend{propo}
Note the slightly stricter bound $\opn{dist}(t,I^c)\leq{}^3\!/_2\tg$ versus the definition of well-preparedness is to leave room for the upcoming convex integration scheme in the next section.

\ssect{Subdividing the time interval}
\ssectcnt*{teor}\ssectcnt*{propo}\ssectcnt*{defi}\ssectcnt*{oss}\ssectcnt*{cor}\ssectcnt*{lemma}
We first introduce a subdivision of the time interval $[0,1]$. Then on each subinterval $[t_i,t_{i+1}]$ we solve a generalized fractional Navier-Stokes system and obtain a solution $v_i$ so that $u+v_i$ is an exact solution of the fractional Navier-Stokes equations \eqref{eq:FNS} on $[t_i,t_{i+1}]$.
\nipar Let $\tg>0$ be a small length scale to be fixed at the end of this section, and define:
\[t_i\coloneq i\tg^\eg\qquad0\leq u\leq\lfloor\tg^{-\eg}\rfloor.\]
Without loss of generality, we assume $\tg^{-\eg}$ is always an integer so that the time interval $[0,1]$ is perfectly divided.
\nipar For $0\leq i\leq\tg^{-\eg}-1$, let $v_i:[t_i,t_{i+1}]\x\T^d\to\R^d$ and $q_u:[t_i,t_{i+1}]\x\T^d\to\R$ be the solution of the following system:
\[\begin{sistema}
\pd_tv_i+\xfrl v_i+\opn{div}(v_i\otimes v_i)+\opn{div}(v_i\otimes u+u\otimes v_i)+\grad q_i=-\opn{div}R \\
\opn{div}v_i=0 \\
v_i(t_i)=0
\end{sistema}. \stag{DiffFNSR}\]
Since the initial data for $v_i$ is zero and $u$ and $R$ are smooth on $[0,1]\x\T^d$, thanks to the general local well-posedness theory for the fractional Navier-Stokes equations, for all sufficiently small $\tg>0$ we may solve the above system on intervals $t\in[t_i,t_{i+1}]$ to obtain a unique smooth solution $v_i$.
\nipar We shall focus on estimating each $v_i$ on the associated interval $[t_i,t_{i+1}]$. The solution $v_i$ serves as an ``accumulator'' of the stress error on $[t_i,t_{i+1}]$, and it will provide the major contribution to the new stress error $\lbar R$ once we use a gluing procedure.
\nipar Recall that $\s{R}:\Cinf(\T^d,\R^d)\to\Cinf(\T^d,\s{C}^{d\x d}_0)$ is an inverse divergence operator on $\T^d$ defined in \cite[Appendix B]{CL}. The below result quantifies the size of the corrector $v_i$ in relation to the time scale $\tg$ and the forcing $-\opn{div}R$.
\xbegin{propo}[][thm:propo:Estvizi]
Let $d\geq2$ and $(u,R)$ be a smooth solution of \eqref{eq:FNSR}. There exists a universal constant $C_r$ depending on $1<r<\8$ so that the following holds.
\nipar For any $\dg>0$, if $\tg>0$ is sufficiently small, then the unique smooth solution $v_i$ to \eqref{eq:DiffFNSR} on $[t_i,t_{i+1}]$ satisfies
\begin{align*}
\norm{v_i}_{L^\8_tH^d_x([t_i,t_{i+1}]\x\T^d)}\leq{}&\dg \\
\norm{\s{R}v_i}_{L^\8_tL^r_x([t_i,t_{i+1}]\x\T^d)}\leq{}&C_i\xints{t_i}{t_{i+1}}\2\norm{R(t)}_{L^r}\diff t+C_u\dg\tg^\eg,
\end{align*}
where $C_u$ is a sufficiently large constant depending on $u$ but not $\dg$ or $\tg$.
\xend{propo}
\begin{qeddim}
The first estimate follows from a simple $H^d$ energy estimate combined with Grönwall's inequality.
\nipar To prove the second one, assume $\dg>0$ is sufficiently small, and denote $v\coloneq v_i,z\coloneq\s{R}v$ for brevity. Note that \eqref{eq:DiffFNSR} preserves the zero-mean condition.
\nipar Let $\P$ be the Leray projection onto divergence-free vector fields on $\T^d$. By projecting \eqref{eq:DiffFNSR} with $\P$ and applying $\s{R}$, we find that the evolution of $z$ is governed by
\[\pd_tz+\xfrl z=-\s{R}\P\opn{div}(v\otimes v+u\otimes v+v\otimes u)-\s{R}\P\opn{div}R.\]
In order to obtain an estimate in $L^r$, we multiply both sides by $z|z|^{r-2}$:
\[z\pd_tz|z|^{r-2}+z\xfrl z|z|^{r-2}=-\s{R}\P\opn{div}(v\otimes v+v\otimes u+u\otimes v)z|z|^{r-2}-\s{R}\P\opn{div}Rz|z|^{r-2}.\]
Integrating on $\T^d$, we obtain
\begin{align*}
\xfr1r\pd_t\norm{z}_{L^r}^r={}&-\xints{\T^d}{}\2(\xfrl z)(z|z|^{r-2})\diff x+\xints{\T^d}{}\2\s{R}\P\opn{div}Rz|z|^{r-2}\diff x \\
&{}-\xints{\T^d}{}\s{R}\P\opn{div}(v\otimes v+v\otimes u+u\otimes v)z|z|^{r-2}\diff x.
\end{align*}
By applying Hölder's inequality, this implies
\[\norm{z}_{L^r}^{r-1}\pd_t\norm{z}_{L^r}\leq\pa{\norm{\xfrl z}_{L^r}+\norm{\s{R}\P\opn{div}(v\otimes v+v\otimes u+u\otimes v}_{L^r}+\norm{\s{R}\P\opn{div}R}_{L^r}}\norm{z|z|^{r-2}}_{L^{\xfr{r}{r-1}}}.\]
Since $\|z|z|^{r-2}\|_{L^{\xfr{r}{r-1}}}=\|z\|_{L^r}^{r-1}$ and $\s{R}\P\opn{div}$ is a Canderón-Zygmund operator, we conclude
\[\pd_t\norm{z}_{L^r}\lsim\norm{\xfrl z}_{L^r}+\norm{v\otimes v+v\otimes u+u\otimes v}_{L^r}+\norm{R}_{L^r}.\]
Integrating over $[t_i,t]$, we obtain
\begin{align*}
\norm{z(t)}_{L^r}\lsim{}&\norm{z(t_i)}_{L^r}+\xints{t_i}{t}\2\norm{R(s)}_{L^r}+\norm{\xfrl z(s)}_{L^r}\diff s \\
&{}+\xints{t_i}{t}\2\norm{(v\otimes v)(s)}_{L^r}+\norm{(v\otimes u)(s)}_{L^r}+\norm{(u\otimes v)(s)}_{L^r}\diff s.
\end{align*}
$z(t_i)=\s{R}(v(t_i)\!)=0$. The last term is easily bounded by the first one in the desired estimate. All other terms must be estimated by $C_u\dg\tg^\eg$. Since they are integrated over an interval of length at most $\tg^\eg$, the following claim concludes the proof:
\[\norm{v\otimes v}_{L^r}+\norm{v\otimes u}_{L^r}+\norm{u\otimes v}_{L^r}+\norm{\xfrl z}_{L^t}\lsim\dg,\]
for a.e. $s\in[t_i,t_{i+1}]$, where the constant can depend on $u$. Since $u$ is smooth, $\|u\|_{L^\8L^\8}$ is finite, so
\[\norm{v\otimes u}_{L^r}+\norm{u\otimes v}_{L^r}\leq\norm{u}_{L^\8 L^\8}\norm{v}_{L^\8 L^r}\leq\norm{u}_{L^\8 L^\8}\norm{v}_{L^\8H^d}\lsim\dg.\]
Analogously
\[\norm{v\otimes v}_{L^r}\leq\norm{v}_{L^\8L^\8}\norm{v}_{L^\8L^r}\leq\norm{v}_{L^\8H^d}^2\leq\dg^2\leq\dg,\]
provided $\dg\leq1$. Coming to the fractional laplacian term, define
\[I(\EXP)\coloneq\xfr12\lfloor2\EXP\rfloor\qquad D(\EXP)\coloneq\EXP-I(\EXP).\]
We can then write:
\begin{align*}
\norm{\xfrl z}_{L^r}^2\lsim{}&\norm{\xfrl v}_{L^2}^2=\norm{\frl{I(\EXP)}[\frl{D(\EXP)}v]}_{L^2}^2 \\
{}={}&\e[13]\sum_{i_1,\dotsc,i_{I(\EXP)}=1}^d\e\norm{\pd_{i_1}^2\dotso\pd_{i_{I(\EXP)}}^2\frl{D(\EXP)}v}_{L^2}^2 \\
{}\lsim{}&\e[13]\sum_{i_1,\dotsc,i_{I(\EXP)}=1}^d\e\abs{\pd_{i_1}^2\dotso\pd_{i_{I(\EXP)}}^2v}_{\dot H^{2D(\EXP)}}^2\leq\norm{v}_{H^{2\EXP}}^2.
\end{align*}
Here, we have used:
\bi
\item The fact $r\leq2$ and thus, since $\T^d$ is of finite measure, $L^2\hra L^r$, and that $\s{R}:L^r\to L^r$ is bounded;
\item A composition property of fractional laplacians;
\item The definition of integer laplacian;
\item A combination of \cite[Proposition 3.3, p. 14]{FracSob} and \cite[Proposition 3.6, p. 18]{FracSob}.
\ei
By \cite[Proposition 2.1, p. 6]{FracSob}, assuming $2\EXP\leq d$, we have:
\[\norm{\xfrl z}_{L^2}\leq\norm{v}_{H^d},\]
which we have already estimated by $C_u\dg$. This completes the proof. \placeqed*
\end{qeddim}

\ssect{Temporal concentration by sharp gluing}
Since $u+v_i$ is an exact solution of the fractional Navier-Stokes equations \eqref{eq:FNS} on each interval $[t_i,t_{i+1}]$ for $0\leq i\leq\tg^{-\eg}-1$, the next step is a suitable gluing of the $v_i$ so that the glued solution $u+\sum_i\xg_iv_i$ is still an exact solution on a majority of the time interval $[0,1]$, with an error supported on many disjoint sub-intervals.
\nipar We first choose cutoff functions to glue the $v_i$ together. We define $\xg_i\in\Cinf_c(\R)$ to be a smooth cutoff such that, for $1\leq i\leq\tg^{-\eg}-2$,
\[\xg_i=\case{
1 & t_i+\tg\leq t\leq t_{i+1}-\tg \\
0 & t_i+\xfr\tg2\geq t\text{ or }t\geq t_{i+1}-\xfr\tg2
}, \stag{Defxi1}\]
and for $i\in\{0,\tg^{-\eg}-1\}$
\begin{align*}
\xg_0={}&\case{
1 & 0\leq t\leq t_1-\tg \\
0 & t\geq t_1-\xfr\tg2
}, \\
\xg_{\tg^{-\eg}-1}=\case{
1 & t_{\tg^{-\eg}-1}+\tg\leq t\leq1 \\
0 & t\leq t_{\tg^{-\eg}}+\xfr\tg2
}. \stag{Defxi2}
\end{align*}
In other words, we do not cut near the endpoints $t=0,t=1$, and the glued solution $\lbar u$ is an exact solution near $t=0$ and $t=1$. It is worth noting that in the iteration scheme $v_i$ for $i=0,i=\tg^{-\eg}-1$ will be zero after step 1, since it already solves \eqref{eq:FNSR} exactly there, and thus the above properties of $\xg_0,\xg_{\tg^{-\eg}-1}$ are only used once.
\nipar Furthermore, we require the following bounds uniformly in $\tg,i$:
\[\abs{\grad^m\xg_i}\lsim_m\tg^{-m}.\]
Note that for sub-intervals $[t_i,t_{i+1}]$, with $1\leq i\leq\tg^{-\eg}-2$, we cut near both the left and the right endpoint. The left cutoff is to ensure smoothness near $t_i$, since each $v_i$ only has a limited amount of time regularity at $t=t_i$, whereas the right cutoff is where the gluing will take place. With $\xg_i$ in hand, we can simply define the glued solution:
\[\lbar u\coloneq u+\sum_i\xg_iv_i=u+\lbar w.\]
It is clear that $\lbar u:[0,1]\x\T^d\to\R^d$ is spatially mean-free and solenoidal. It remains to show that $\lbar u$ satisfies the properties listed in \kcref{thm:propo:StressConc}.
\nipar Heuristically, $\lbar u$ should be an exact solution with a stress error supported on smaller intervals of size $\tg$. To confirm this claim, we must computer the stress error $\lbar R$ associated with $\lbar u$. Since the $\xg_i$ are disjointly supported, we can compute
\begin{align*}
\pd_t\lbar u+\xfrl\lbar u={}&\opn{div}R+(\pd_t+\xfrl)\sum_i\xg_iv_i \\
&{}+\sum_i\xg_i\opn{div}(u\otimes v_i)+\sum_i\xg_i\opn{div}(v_i\otimes u) \\
&{}+\sum_i\xg_i^2\opn{div}(v_i\otimes v_i).
\end{align*}
Thus, using the fact that $v_i$ solves \eqref{eq:DiffFNSR} on $[t_i,t_{i+1}]$ and $u$ solves \eqref{eq:FNSR} on $[0,1]$, we have:
\begin{align*}
\pd_t\lbar u+\xfrl u+\opn{div}(\lbar u\otimes\lbar u)+\grad p={}&\opn{div}R+\sum_i\pd_t\xg_iv_i+\sum_i(\xg_i^2-\xg_i)\opn{div}(v_i\otimes v_i) \\
&{}+\sum_i\xg_i(\pd_tv_i-\xfrl v_i+\opn{div}(\!(u+v_i)\otimes v_i)+\opn{div}(u_i\otimes u)\!) \\
{}={}&\pa{1-\sum_i\xg_i}\opn{div}R+\sum_i\pd_t\xg_iv_i \\
&{}+\sum_i(\xg_i^2-\xg_i)\opn{div}(v_i\otimes v_i)-\sum_i\xg_i\grad q_i.
\end{align*}
Now let
\[\lbar R\coloneq\pa{1-\sum\xg_i}R+\s{R}\sum_i\pd_t\xg_iv_i+\sum_i(\xg_i^2-\xg_i)v_i\mathbin{\pint\otimes}v_i, \stag{defRbar}\]
where $\pint\otimes$ denotes a traceless tensor produce, i.e. $f\mathbin{\pint\otimes}g=f_ig_j-\xfr1d\dg_{ij}f_kg_k$. Since each $v_i$ has zero spatial mean, $\opn{div}\s{R}v_i=v_i$, and we can conclude that
\[\pd_t\lbar u-\xfrl\lbar u+\opn{div}(\lbar u\otimes\lbar u)+\grad\lbar p=\opn{duv}\lbar R,\]
where the pressure $\lbar p:[0,1]\x\T^d\to\R$ is defined by:
\[\lbar p=p+\sum_i\xg_iq_i-\sum_i(\xg_i^2-\xg)\xfr{|v_i|^2}{d}.\]
The last step is then to show that the new Reynolds stress $\lbar R$ is comparable to the original one in the $L^1_tL^r$ norm. It is clear that $\lbar R$ is much more ``turbulent'' than the original $R$ as its value changes much more drastically due to the sharp cutoffs near the endpoints of each $[t_i,t_{i+1}]$.
\nipar The heuristic is that, if $\tg$ is small enough, then $v_i$ linearly with a rate of order $\opn{div}R$, and therefore gluing the $v_i$ together only counts the input from the stress forcing $\opn{div}R$. More precisely, the leading order term in \eqref{eq:defRbar} is the second term, where $\s{R}v_i$ is proportional to $R$ thanks to \kcref{thm:propo:Estvizi}.
\xbegin{propo}
For any $1<r<\8$, there exists a universal constant $C_r$ depending on $r$ and $\eg$ such that for all sufficiently small $\tg>0$ the glue solution $(\lbar u,\lbar R)$ safistifs
\[\norm{\lbar R}_{L^1_tL^r_x}\leq C_r\norm{R}_{L^1_tL^r_x}.\]
\xend{propo}
The proof of this is left to \cite{CL}, where it is Proposition 3.3 on p. 16.

\ssect{Proof of \texorpdfstring{\kcref{thm:propo:StressConc}}{Propo 3.1 CL}}
We conclude this section with the last step in the proof of \kcref{thm:propo:StressConc}. Since all the estimates have been obtained, we only need to verify that the temporal support of $\lbar w=\sum_i\xg_iv_i$ is contained in $\tilde I$, and that $(\lbar u,\lbar R)$ is well-prepared.
\nipar Note that $(u,R)$ is well-prepared for $\tilde I,\tilde\tg$, and it follows from \eqref{eq:DiffFNSR} that
\[v_i\equiv0\qquad0\leq i\leq\tg^{-\eg}\wedge R|_{[t_i,t_{i+1}]}\equiv0.\]
Hence, if $\tg^\eg=|[t_i,t_{i+1}]|$ is sufficiently smaller than $\tilde\tg$, the definition of well-preparedness of $(u,R)$ implies that
\[\bigcup_i\opn{supp}_t(\xg_iv_i)\sbs\tilde I.\]
Thus we have proved that $\opn{supp}_t\lbar w\sbs\tilde I$.
\nipar Let us now show the well-preparedness of $(\lbar u,\lbar R)$. Define an index set
\[E\coloneq\br{i\in\Z:1\leq i\leq\tg^{-\eg}-1,v_i\not\equiv0},\]
satisfying a trivial estimate:
\[\abs{E}\leq\tg^{-\eg}.\]
The idea is that the concentrated stress $\lbar R$ is supported around each $t_i$ for $i\in E$. Therefore, we can define a set on the time axis
\[I\coloneq\bigcup_{i\in*E}\sq{t_i-\xfr52\tg,t_i+\xfr52\tg},\]
where as before $t_i\coloneq i\tg^\eg$. Note that each interval in $I$ has length $5\tg$ and the total number of intervals is at most $\tg^{-\eg}$, consistent with the well-preparedness, and $0,1\nin I$ due to \eqref{eq:Defxi2}.
\nipar Now take any $t\in[0,1]$ such that $\opn{dist}(t,I^c)\leq\xfr32\tg$. Then by \eqref{eq:Defxi1} and \eqref{eq:Defxi2}
\[\sum_i\xg_i(t)=1.\]
Moreover,  $\pd_t\xg_i(t)=0$ and $\xg_i(t)\in\{0,1\}$ for any $i$. Consequently,
\[\lbar R(t)=\pa{1-\sum_i\xg_i}R+\s{R}\sum_i\pd_t\xg_iv_i+\sum_i(\xg_i^2-\xg_i)v_i\mathbin{\pint\otimes}v_i=0,\]
for every $t$ such that $\opn{dist}(t,I^c)\leq\xfr32\tg$. In particular, $(\lbar u,\lbar R)$ is also well-prepared, which concludes the proof.
\xbegin{oss}[First relations]
You may have noted that we obtained an estimate for $v$ in $L^\8_tH^d_x$, but our meta-theorem requires estimates in $L^p_tL^\8_x\cap L^s_tW^{\gg,q}_x$. Now, $L^\8_tH^d_x\hra L^p_tL^\8_x$ for any $p\geq1$ since we are on $\T^d$, but the embedding $L^\8_tH^d_x\hra L^s_tW^{\gg,q}_x$ is not guaranteed, since the embedding $H^d_x\hra W^{\gg,q}_x$ is not. It is clear that, if $\gg\leq d$ and $q\leq 2$, or if $\gg\leq\xfr d2$, the embedding holds. By fractional Sobolev embedding, this will hold if:
\[q\leq\xfr{2d}{2\gg-d}\quad\text{or}\quad2\gg\leq d.\]
We can abbreviate this by writing:
\[\gg\leq d\wedge q\leq\xfr{2d}{(2\gg-d)^+}, \stag{QG}\]
where ``<1/0'' is understood to mean ``arbitrary''.
\end{oss}

\sect{Convex integration in space-time}\label{Conc}
\sectcnt*{teor}\sectcnt*{propo}\sectcnt*{defi}\sectcnt*{oss}\sectcnt*{cor}\sectcnt*{lemma}
In this section, we will use a convex integration scheme to reduce the size of the Reynolds stress. The goal is to design a suitable velocity perturbation $w$ to the glued solution $(\lbar u,\lbar R)$ so that $u_1\coloneq\lbar u+w$ solves \eqref{eq:FNSR} with a much smaller Reynolds stress $R_1$.
\nipar The main goal of the current and the following section is summarized in the following proposition.
\xbegin{propo}[Main Perturbation Step][thm:propo:MPS]
There exists a universal constant $M>0$ such that for any $p,q,s,\gg,\EXP$ in the ranges of the main theorem and $d\geq3$ there exists $r=r(p,q,s,d)>1$ such that the following holds.
\nipar Let $\dg>0$ and $(\lbar u,\lbar R)$ be a well-prepared smooth solution of \eqref{eq:FNSR}. Then there exists another well-prepared smooth solution $(u_1,R_1)$ of \eqref{eq:FNSR} for the same set $I$ and time scale $\tg$ such that:
\[\norm{R_1}_{L^1_tL^r_x}\leq\dg.\]
Moreover, the velocity perturbation $w\coloneq u_1-u$ satisfies:
\begin{align}
\opn{supp}w\sbs{}&I\x\T^d \label{eq:SuppPert} \\
\norm{w}_{L^2_tL^2_x}\leq{}&M\norm{\lbar R}_{L^1_tL^1_x} \label{MPS:L2L2} \\
\norm{w}_{L^p_tL^\8_x}+\norm w_{L^s_tW^{\gg,q}_x}\leq{}&\dg. \label{MPS:LpLq+LsW1q}
\end{align}
\xend{propo}
In the remainder of this section, we will construct the velocity perturbation $w$ and define its associated Reynolds stress $R_1$ and pressure $p_1$. The well-preparedness of $(u_1,R_1)$ will be an easy consequence of the definition of $w$, whereas all the estimates will be proven in the next section.

\ssect{Stationary Mikado flows for the convex integration}
\ssectcnt*{teor}\ssectcnt*{propo}\ssectcnt*{defi}\ssectcnt*{oss}\ssectcnt*{cor}\ssectcnt*{lemma}
The main building blocks of the convex integration scheme are the Mikado flows $W_k:\T^d\to\R^d$ introduced in \cite{DSz}. In other contexts, they are called concentrated Mikado flows, or Mikado flows with concentration since they are supported on periodic cylinders with a small radius. Here for brevity, we refer to them as the Mikado flows.
\nipar We start with a geometric lemma that dates back to the work of Nash and is the key reason to use Mikado flows in convex integration schemes. $\s{S}^{d\x d}_+$ is the set of positive definite $d\x d$ matrices and $e_k\coloneq\xfr{k}{|k|}$ for all $k\in\Z^d$.
\xbegin{lemma}[Decomposition lemma][thm:lemma:DecLemma]
For any compact subset $\s{N}\sbs\s{S}^{d\x d}_+$, there exist a finite set $\Lg\sbs\Z^d$ and smooth functions $\Gg_k\in\Cinf(\s{N};\R)$ for any $k\in\Lg$ such that:
\[R=\sum_{k\in*\Lg}\Gg_k^2(R)e_k\otimes e_k\qquad\VA R\in\s{N}.\]
\xend{lemma}
We apply this lemma for $\s{N}=B_{\xfr12}(\opn{Id})$ (the metric ball of radius $\xfr12$ around the identity $\opn{Id}$ in the space $\s{S}^{d\x d}_+$) to obtain smooth functions $\Gg_k$ for $k\in\Lg\sbs\Z^d$. Throughout this note, the direction set $\Lg$ is fixed, and we construct the Mikado flows as follows.
\nipar We choose points $p_k\in(0,1)^d$ such that $p_k\neq p_{-k}$ if both $k,-k\in\Lg$. For each $k\in\Lg$, we denote by $\ell_k\sbs\T^d$ the periodic line passing through $p_k$ in direction $k$, namely:
\[\ell_k\coloneq\br{tk+p_k\in\T^d:t\in\R}.\]
Since $\Lg$ is a finite lattice set and we identify $\T^d$ with a periodic box $[0,1]^d$, there exists a geometric constant $C_\Lg\in\N$ depending on the set $\Lg$ such that:
\[\abs{\ell_k\cap\ell_{k'}}\leq C_\Lg\qquad\VA k,k'\in\Lg,\]
where we note that $\ell_k\cap\ell_{-k}=\0$ due to $p_k\neq p_{-k}$.
Let $\mg>0$ be the spatial concentration parameter whose value will be fixed in the end of the proof. Let $\fg,\yg\in\Cinf_c([\xfr12,1])$ be such that, if we define $\yg_k,\fg_k:\T^d\to\R$ by:
\begin{align*}
\yg_k\coloneq{}&\mg^{\xfr{d-1}{2}}\yg(\mg\opn{dist}(\ell_k,x)\!) \\
\fg_k\coloneq{}&\mg^{\xfr{d-1}{2}-2}\fg(\mg\opn{dist}(\ell_k,x)\!), \\
\intertext{then:}
\Dg\fg_k={}&\yg_k\quad\text{on }\T^d\quad\text{and}\quad\fint_{\T^d}\yg_k^2\diff x=1. \\
\intertext{Note that:}
\opn{supp}\yg_k\cap\opn{supp}\yg_{k'}\sbs{}&\br{x\in\T^d:\opn{dist}(x,\ell_k\cap\ell_{k'})\leq M_\Lg\mg^{-1}},
\end{align*}
for a sufficiently large constant $M_\Lg$ depending on $\Lg$.
\nipar Finally, the stationary Mikado flows $W_k:\T^d\to\R^d$ are defined by:
\[W_k\coloneq\yg_ke_k,\]
where the constant vector $e_k=\xfr{k}{|k|}$. Using the gradient field $\grad\fg_k$, we may write $W_k$ as a divergence of a skew-symmetric tensor $\Wg_k\in\Cinf_0(\T^d,\R^{d\x d})$:
\[\Wg_k\coloneq e_k\otimes\grad\fg_k-\grad\fg_k\otimes e_k.\]
Indeed, $\Wg_k$ is a skew-symmetric tensor by definition, and by a direct computation:
\[\opn{div}\Wg_k=\opn{div}(\grad\fg_k)e_k-(e_k\per\grad)\grad\fg_k=\Dg\fg_ke_k-0=W_k.\]
We summarize the main properties of the Mikado flows $W_k$ in the following theorem.
\xbegin{teor}[Properties of the Mikado flows][thm:teor:PropMikFlows]
Let $d\geq2$ be the dimension. The stationary Mikado flows $W_k:\T^d\to\R^d$ satisfy the following.
\ben[label={\upshape(\arabic*)}]
\item Each $W_k\in\Cinf_0(\T^d)$ is divergence-free, satisfies:
\[W_k=\opn{div}\Wg_k,\]
and solves the pressureless Euler equations:
\[\opn{div}(W_k\otimes W_k)=0;\]
\item For any $1\leq p\leq\8$, the following estimates hold uniformly in $\mg$:
\[\mg^{-m}\norm{\grad^mW_k}_{L^p(\T^d)}\lsim_m\mg^{\xfr{d-1}{2}-\xfr{d-1}{p}}\qquad\mg^{-m}\norm{\grad^m\Wg_k}_{L^p(\T^d)}\lsim\mg^{-1+\xfr{d-1}{2}-\xfr{d-1}{p}};\]
\item For any $k\in\Lg$, there holds:
\[\fint_{\T^d}\2W_k\otimes W_k=e_k\otimes e_k,\]
and for any $1\leq p\leq\8$:
\[\norm{W_k\otimes W_{k'}}_{L^p(\T^d)}\lsim\mg^{d-1-\xfr dp}\qquad k\neq k'.\]
\een
\xend{teor}
This is \cite[Theorem 4.3]{CL}. Claim (3) can be found proven there. The other two claims are straightforwardly deduced from the above definitions and estimates.
\xbegin{cor}[Mikado flows and fractional Sobolev norms]
By interpolation, claim (2) above implies its fractional version:
\[\norm{W_k}_{W^{s,p}(\T^d)}\lsim\mg^{s+\xfr{d-1}{2}-\xfr{d-1}{p}} \qquad \norm{\Wg_k}_{W^{s,p}(\T^d)}\lsim\mg^{s-1+\xfr{d-1}{2}-\xfr{d-1}{p}}.\]
\end{cor}

\ssect{Implementation of temporal concentration}
Since Mikado flows are stationary, the velocity perturbation will be homogeneous in time if we simply use Decomposition \kcref{thm:lemma:DecLemma} to define the coefficients. To obtain $L^p_tL^\8_x\cap L^s_tW^{\gg,q}_x$ estimates, it is necessary to introduce temporal concentration in the perturbation.
\nipar To this end, we choose temporal functions $g_\kg,h_\kg$ to oscillate the building blocks $W_k$ intermittently in time. Specifically, $g_\kg$ will be used to oscillate $W_k$ so that the space-time cascade balances the low temporal frequency part of the old stress error $\lbar R$, whereas $h_\kg$ is used to define a temporal corrector whose time derivative will further balance the high temporal frequency part of the old stress error $\lbar R$.
\nipar Let $g\in\Cinf_c([0,1])$ be such that:
\[\xints01g^2(t)\diff t=1.\]
To add in temporal concentration, let $\kg>0$ be a large constant whose value will be specified later and define $g_\kg:[0,1]\to\R$ as the 1-periodic extension of $\kg^{\xfr12}g(\kg t)$ so that:
\[\norm{g_\kg}_{L^p([0,1])}\lsim\kg^{\xfr12-\xfr1p}\qquad\VA p\in[1,\8].\]
The value of $\kg$ will be specified later and the function $g_\kg$ will be used in the definition of the velocity perturbation. As we will see in \kcref{thm:lemma:OscReynPress}, the nonlinear term can only balance a portion of the stress $\lbar R$ and there is a leftover term which is of high temporal frequency. This motivates us to consider the following temporal corrector.
\nipar Let $h_\kg:[0,1]\to\R$ be defined by:
\[h_\kg(t)=-t+\xints0tg_\kg^2(s)\diff s.\]
For $\kg\in\N$, in view of the zero-mean condition for $g_\kg^2(t)-1$, the function $h_\kg:[0,1]\to\R$ is well-defined and periodic, and we have:
\[\norm{h_\kg}_{L^\8([0,1])}\leq1,\]
uniformly in $\kg$.
\nipar We remark that for any $\ng\in\N$, the periodically rescaled function $g_\kg(\ng\per):[0,1]\to\R$ also verifies the bound:
\[\norm{g_\kg(\ng\per)}_{L^p([0,1])}\lsim\kg^{\xfr12-\xfr1p}\qquad\VA p\in[1,\8].\]
Moreover, we have the identity:
\[\pd_t(\ng^{-1}h_\kg(\ng t)\!)=g_\kg^2-1,\]
which will imply the smallness of the corrector, cfr. the definition of $w\stp t$ below.

\ssect{Space-time cutoffs}
Before introducing the velocity perturbation, we need to define two important cutoff functions, one to ensure \kcref{thm:lemma:DecLemma} applies and the other to ensure the well-preparedness of the new solution $(u_1,R_1)$.
\nipar Since \kcref{thm:lemma:DecLemma} is stated for a fixed compact set in $\s{S}^{d\x d}_+$, we need to introduce a cutoff for the stress $\lbar R$. Let $\xg:\R^{d\x d}\to\R^+$ be a positive smooth function such that $\xg$ is monotonically increasing with respect to $|x|$ and:
\[\xg(x)=\case{
1 & 0\leq|x|\leq1 \\
|x| & |x|\geq2
}.\]
With this cutoff $\xg$, we may define a divisor for the stress $\lbar R$ so that \kcref{thm:lemma:DecLemma} applies. Indeed, define $\rg\in\Cinf([0,1]\x\T^d)$ by:
\[\rg=2\xg(\lbar R).\]
Then immediately:
\[\opn{Id}-\xfr{\lbar R}{\rg}\in B_{\xfr12}(\opn{Id})\qquad\VA(t,x)\in[0,1]\x\T^d,\]
which means we can use $\opn{Id}-\xfr{\lbar R}{\rg}$ as the argument in the smooth functions $\Gg_k$ given by \kcref{thm:lemma:DecLemma}.
\nipar Next, we need another cutoff to take care of the well-preparedness of the new solution $(u_1,R_1)$ as the perturbation has to be supported within $I$. Let $\qg\in\Cinf_c(\R)$ be a smooth temporal cutoff function such that:
\[\qg(t)=\case{
1 & \opn{dist}(t,I^C)\geq\xfr32\tg \\
0 & \opn{dist}(t,I^C)\leq\tg
},\]
where $I\sbs[0,1]$ and $\tg>0$ are the parameters of the well-preparedness of the glued solution $(\ubar u,\ubar R)$, thus required to satisfy $I\sbse\tilde I,0,1\nin I,\tg<2^{-1}\tilde\tg$, where $\tilde I,\tilde\tg$ are the parameters of the well-preparedness of the unglued solution $(u,R)$. Note that this cutoff ensures that the new solution will still be well-prepared.

\ssect{The velocity perturbation}
We recall the four parameters for the perturbation we have defined so far, and add a fourth one:
\ben[label=(\arabic*)]
\item Temporal oscillation $\ng\in\N$;
\item Temporal concentration $\kg>0$;
\item Spatial concentration $\mg\in\N$;
\item Spatial oscillation $\sg\in\N$, which we introduce here, serves the purpose of making the $W_k$ oscillate faster.
\een
The requirements $\ng,\sg,\mg\in\N$ are for the sake of periodicity. We will now define and estimate the new Reynolds and the velocity perturbation in terms of these parameters, and obtain a series of relations that yield the bounds we desire, and then we will see for what ranges of $s,p,q,\gg,\EXP$ we can choose the parameters so that the relations are satisfied.
\nipar With all the ingredients in hand, we are ready to define the velocity perturbation. In summary, the velocity perturbation $w:[0,1]\x\T^d\to\R^d$ consists of three parts:
\[w=w\stp p+w\stp c+w\stp t.\]
The principal part of the perturbation $w\stp p$ consists of super-positions of the building blocks $W_k$ oscillating with period $\sg^{-1}$ on $\R^d$ and period $\ng^{-1}$ on $[0,1]$:
\[w\stp p(t,x)\coloneq\sum_{k\in*\Lg}a_k(t,x)W_l(\sg x),\]
where the amplitude function $a_k:[0,1]\x\T^d\to\R$ is given by:
\[a_k\coloneq\qg g_\kg(\ng t)\rg^{\xfr12}\Gg_k\pa{\opn{Id}-\xfr{\lbar R}{\rg}}.\]
Note that the above $w\stp p$ is not divergence-free. To fix this, we introduce a divergence-free corrector using the tensor potential $\Wg_k$:
\[w\stp c(t,x)\coloneq\sg^{-1}\per\!\sum_{k\in*\Lg}\grad a_k(t,x)\mathbin:\Wg_k(\sg x).\]
Indeed, we can rewrite:
\begin{align*}
w\stp p+w\stp c={}&\sg^{-1}\per\!\sum a_k(t,x)\opn{div}\Wg_k(\sg x)+\sg^{-1}\per\!\sum\grad a_k(t,x)\mathbin:\Wg_k(\sg x)={} \\
{}={}&\sg^{-1}\opn{div}\sum a_k(t,x)\Wg_k(\sg x),
\end{align*}
where each $a_k\Wg_k$ is skew-symmetric and hence double-divergence-free, making it so that $\opn{div}(w\stp p+w\stp c)=0$.
\nipar Finally, we define a temporal corrector to balance the high temporal frequency part of the interaction. This Ansatz was first introduced in \cite{BV} and also used in \cite{BCV}. The heart of the argument is to ensure that:
\begin{align*}
\pd_tw\stp t+\opn{div}(w\stp p\otimes w\stp p)={}&\text{Pressure gradient} \\
&{}+{}\text{Terms with high spacial frequencies} \\
&{}+{}\text{Lower order terms}.
\end{align*}
However, the key difference between \cite{BV,BCV} and the current scheme is that here the smallness of the corrector is free and it does not require much temporal oscillation, which is the reason we must use stationary spatial building blocks.
\nipar Specifically, the temporal corrector $w\stp t$ is defined as:
\[w\stp t\coloneq\ng^{-1}h_\kg(\ng t)(\opn{div}(\lbar R)-\grad\Dg^{-1}\opn{div}\opn{div}(\lbar R)\!),\]
where we note that $\Dg^{-1}$ is well-defined on $\T^d$ since $\opn{div}\opn{div}(\lbar R)=\pd_i\pd_j\lbar R_{ij}$ has zero spatial mean.
\nipar It is easy to check $\opn{supp}_tw\stp t\sbs\opn{supp}_t\lbar R$ and $w\stp t$ is divergence-free. Indeed:
\[\opn{div}w\stp t=\ng^{-1}h_\kg(\ng t)(\pd_i\pd_j\lbar R_{ij}-\pd_k\pd_k\Dg^{-1}\pd_i\pd_j\lbar R_{ij})=0.\]
In the lemma below, we show that the leading order interaction of the principal part $w\stp p$ is able to balance the low temporal frequency part of the stress error $\lbar R$, which motivates the choice of the corrector $w\stp t$.
\xbegin{lemma}[Properties of the coefficients][thm:lemma:Propak]
The coefficients $a_k$ satisfy:
\[a_k=0\quad\text{if}\quad\opn{dist}(t,I^C)\leq\tg,\]
and:
\[\sum_{k\in*\Lg}a_k^2\fint_{\T^d}W_k\otimes W_k\diff x=\qg^2g_\kg^2(\ng t)\rg\opn{Id}-g_\kg^2(\ng t)\lbar R.\]
\xend{lemma}
The proof of this can be found in \cite{CL}, where it is Lemma 4.4.

\ssect{The new Reynolds stress}
In this subsection, our goal is to design a suitable stress tensor $R_1:[0,1]\x\T^d\to\s{S}_0^{d\x d}$ such that the pair $(u_1,R_1)$ is a smooth solution of \eqref{eq:FNSR} for a suitable smooth pressure $p_1$.
\nipar We first compute the nonlinear term and isolate nonlocal interactions:
\[\opn{div}(w\stp p\otimes w\stp p+\lbar R)=\pon{div}\sq{\sum_ka_k^2W_k(\sg x)\otimes W_k(\sg x)}+\opn{div}R_{far},\]
where $R_{far}$ denotes the nonlocal interactions between Mikado flows of different directions:
\[R_{far}=\sum_{k\neq k'}a_ka_{k'}W_k(\sg x)\otimes W_{k'}(\sg x).\]
And then we proceed to examine the first term in the above decomposition, for which by \kcref{thm:lemma:Propak} we have:
\begin{multline*}
\opn{div}\sq{\sum_ka_k^2W_k(\sg x)\otimes W_k(\sg x)}={} \\
{}=\opn{div}\sq{\sum_ka_k^2\sq{\pa{W_k(\sg x)\otimes W_k(\sg x)-\fint W_k\otimes W_k}+\fint W_k\otimes W_k}+\lbar R}={} \\
{}=\opn{div}\sum_ka_k^2\pa{W_k(\sg x)\otimes W_k(\sg x)-\fint W_k\otimes W_k}+\grad(\qg^2g_\kg^2\rg)+(1-g_\kg^2(\ng t)\!)\opn{div}\lbar R.
\end{multline*}
Finally, using the product rule, we compute the divergence term as:
\[\opn{div}\sum_ka_k^2\pa{W_k(\sg x)\otimes W_k(\sg x)-\fint W_k\otimes W_k}=\sum_k\grad(a_k^2)\per\pa{W_k(\sg x)\otimes W_k(\sg x)-\fint W_k\otimes W_k}.\]
Typical in the convex integration, we can gain a factor of $\sg^{-1}$ in the above equality by inverting the divergence. To this end, let us use the bilinear anti-divergence operator $\s{B}$ from \kcref{thm:teor:BilinAntidivEq}. Since the above equality has zero spatial mean, by \eqref{eq:Antidiv2} it is equal to $\opn{div}R_{osc,x}$ where:
\[R_{osc,x}=\sum_k\s{B}\pa{\grad(a_k^2),W_k(\sg x)\otimes W_k(\sg x)-\fint W_k\otimes W_k}.\]
Combining all the above, we have:
\[\opn{div}(w\stp p\otimes w\stp p+\lbar R)=\pon{div}R_{osc,x}+\pon{div}R_{far}+\grad(\qg^2g_\kg^2\rg)+(1-g_\kg^2(\ng t)\!)\pon{div}\lbar R.\]
In view of the above computations, we define a temporal oscillation error:
\[R_{osc,t}\coloneq\ng^{-1}h_\kg(\ng t)\s{R}\opn{div}(\pd_t\lbar R),\]
so that the following decomposition holds.
\xbegin{lemma}[Oscillation Reynolds and pressure][thm:lemma:OscReynPress]
Let the space-time oscillation error $R_{osc}$ be:
\[R_{osc}\coloneq R_{osc,x}+R_{osc,t}+R_{far}.\]
Then:
\[\pd_tw\stp t+\opn{div}(w\stp p\otimes w\stp p+\lbar R)+\grad P=\pon{div}R_{osc},\]
where the pressure term $P$ is defined by:
\[P\coloneq\qg^2g_\kg^2(\ng t)\rg-\ng^{-1}\Dg^{-1}\pon{div}\opn{div}\pd_t(\lbar Rh_\kg(\ng t)\!).\]
\xend{lemma}
The proof of this is a simple calculation which is left to \cite{CL}, where it is Lemma 4.5.
\nipar Finally, we can define the correction error and the linear error as usual:
\begin{align*}
R_{cor}\coloneq{}&\s{R}(\opn{div}(w\stp c+w\stp t)\otimes w+w\stp p\otimes(w\stp c+w\stp t)\!) \\
R_{lin}\coloneq{}&\s{R}(\pd_t(w\stp p+w\stp c)-\Dg\lbar u+\lbar u\otimes w+w\otimes\lbar u).
\end{align*}
To conclude, we summarize the main results in this section below.
\xbegin{lemma}[The new Reynolds stress][thm:lemma:NewReyn]
Define the new Reynolds stress and pressure by:
\begin{align*}
R_1\coloneq{}&R_{lin}+R_{cor}+R_{osc} \\
p_1\coloneq{}&\lbar p-P.
\end{align*}
Then $(u_1,R_1)$ is a well-prepared solution to NSR and the velocity perturbation $w\coloneq u_1-\lbar u$ satisfies $\opn{supp}w\sbs I\x\T^d$.
\xend{lemma}
The proof of this is also left to \cite{CL}, where it is Lemma 4.6.

\sect{Proof of the iteration proposition}\label{Est}
\sectcnt*{teor}\sectcnt*{propo}\sectcnt*{defi}\sectcnt*{oss}\sectcnt*{cor}\sectcnt*{lemma}
In this section we will show that the velocity perturbation $w$ and the new Reynolds stress $R_1$ derived in the previous section satisfy the properties claimed in \kcref{thm:propo:MPS}.
\nipar As a general note, we use a constant $C_u$ for dependency on the previous solution $(\lbar u,\lbar R)$ throughout this section.
\nipar We now work out what relations between $\mg,\sg,\ng,\kg$ we need for our estimates to hold, and then we see how we can achieve them.

\ssect{Estimates on the velocity perturbation}
We first estimate the coefficients $a_k$ of the perturbation $w$. The following is \cite[Lemma 5.2]{CL}.
\xbegin{lemma}[Estimates on the perturbation coefficients][thm:lemma:Estak]
The coefficients $a_k$ are smooth on $[0,1]\x\T^d$ and:
\[\norm{\pd_t^n\grad^ma_k}_{L^p_tL^\8_x}\leq C_{u,m,n}(\ng\kg)^n\kg^{\xfr12-\xfr1p}\qquad p\in[1,\8].\]
In addition, the following bound holds for all times $t\in[0,1]$:
\[\norm{a_k(t)}_{L^2(\T^d)}\lsim\qg(t)g_\kg(\ng t)\pa{\xints{\T^d}{}\2\rg(t,x)\diff x}^{\xfr12}.\]
\xend{lemma}
With these estimates of $a_k$ in hand, we start estimating the velocity perturbation. As expected, the principal part $w\stp p$ is the largest among all parts in $w$. The following adapts \cite[Proposition 5.3]{CL}.
\xbegin{propo}[Estimates on the principal part][thm:propo:Estwp]
Assume that
\begin{align*}
\sg\lsim{}&\ng^2 \stag{Rel0} \\
\kg^{\xfr12-\xfr1p}\mg^{\xfr{d-1}{2}}\lsim\lg^{-\hg} \stag{Rel1} \\
\sg^\gg\kg^{\xfr12-\xfr2s}\mg^{\gg+\xfr{d-1}{2}-\xfr{d-1}{q}}\lsim\lg^{-\hg}. \stag{Rel2}
\end{align*}
The principal part $w\stp p$ satisfies:
\[\norm{w\stp p}_{L^2_tL^2_x}\lsim\norm{\lbar R}^{\xfr12}_{L^1_tL^1_x}+C_u\sg^{-\xfr12}\qquad\norm{w\stp p}_{L^p_tL^\8_x}+\norm{w\stp p}_{L^s_tW^{1,q}}\leq C_u\lg^{-\hg}.\]
In particular, for sufficiently large $\lg,\sg$:
\[\norm{w\stp p}_{L^2_tL^2_x}\lsim\norm{\lbar R}^{\xfr12}_{L^1_tL^1_x}\qquad\norm{w\stp p}_{L^p_tL^\8_x}+\norm{w\stp p}_{L^s_tW^{1,q}_x}\leq\xfr\dg4.\]
\xend{propo}
\clearpage
\begin{qeddim*}[Sketch]
\centered*{$L^2_tL^2_x$ estimate}
The calculations in \cite{CL} lead to the following estimate:
\[\norm{w\stp p}_{L^2_tL^2_x}\lsim\norm{\lbar R}^{\xfr12}_{L^1_tL^1_x}+C_u\ng^{-1}+C_u\sg^{-\xfr12}.\]
The desired estimate thus follows from $\ng^{-1}\lsim\sg^{-\xfr12}$, which is equivalent to \eqref{eq:Rel0}.
\centered*{$L^p_tL^\8_x$ estimate}
The calculations in \cite{CL} lead to the following bound:
\[\norm{w\stp p}_{L^p_tL^\8_x}\lsim_u\kg^{\xfr12-\xfr1p}\mg^{\xfr{d-1}{2}}.\]
The desired estimate thus follows from \eqref{eq:Rel1}.
\centered*{$L^s_tW^{\gg,q}_x$ estimate}
We take the $W^{\gg,q}_x$ norm to obtain:
\[\norm{w\stp p(t)}_{W^{\gg,q}(\T^d)}\lsim\sum\norm{a_k(t)}_{\s{C}^\gg}\norm{W_k(\sg\per)}_{W^{\gg,q}}.\]
Taking the $L^s$ norm in time, by \kcref{thm:teor:PropMikFlows} and \kcref{thm:lemma:Estak} we have that:
\begin{align*}
\norm{w\stp p}_{L^s_tW^{\gg,q}_x}\lsim{}&\sum_k\norm{a_k}_{L^s_t\s{C}^\gg_x}\norm{W_k(\sg\per)}_{W^{\gg,q}(\T^d)}\lsim\sg^\gg\mg^{\gg+\xfr{d-1}{2}-\xfr{d-1}{q}}\per\e\sum_k\norm{a_k}_{L^s_t\s{C}^\gg_x}\lsim{} \\
{}\lsim{}&\sg^\gg\kg^{\xfr12-\xfr1s}\mg^{\gg+\xfr{d-1}{2}-\xfr{d-1}{q}}.
\end{align*}
The desired bound then follows from \eqref{eq:Rel3}.
\end{qeddim*}
Next, we estimate the corrector $w\stp c$, which is expected to be much smaller than $w\stp p$ due to the derivative gains from both the fast oscillation $\sg$ and the tensor potential $\Wg_k$. The following adapts \cite[Proposition 5.4]{CL}.
\xbegin{propo}[Estimates on the incompressibility corrector][thm:propo:Estwc]
Assume that:
\begin{align*}
\sg^{-1}\mg^{\xfr{d-1}{2}-1}\kg^{\xfr12-\xfr1p}\lsim{}&\lg^{-\hg} \stag{Rel3} \\
\sg^{-1}\mg^{\xfr{d-1}{2}-\xfr{d-1}{2}\xfr{2-r}{r}-1}=\sg^{-1}\mg^{d-2-\xfr{d-1}{r}}\lsim{}&\lg^{-\hg}. \stag{Rel3b} \\
\sg^{\gg-1}\kg^{\xfr12-\xfr1s}\mg^{\gg-1+\xfr{d-1}{2}-\xfr{d-1}{q}}\lsim{}&\lg^{-\hg}, \stag{Rel4}
\end{align*}
The divergence-free corrector $w\stp c$ then satisfies:
\[\norm{w\stp c}_{L^2_tL^{\xfr{2r}{2-r}}_x}\leq C_u\lg^{-\hg}\qquad\qquad\norm{w\stp c}_{L^p_tL^\8_x}+\norm{w\stp c}_{L^s_tW^{\gg,q}_x}\leq C_u\lg^{-\hg}.\]
In particular, for sufficiently large $\lg$:
\[\norm{w\stp c}_{L^2_tL^{\xfr{2r}{2-r}}_x}\leq\norm{\lbar R}_{L^1_tL^1_x}^{\xfr12}\qquad\qquad\norm{w\stp c}_{L^p_tL^\8_x}+\norm{w\stp c}_{L^s_tW^{\gg,q}_x}\leq\xfr\dg4.\]
\end{propo}
\xbegin{oss}
An estimate in $L^2_tL^\8_x$ was proved at this point in \cite{CL}. However, that required a relation which severely limited the bounds on $s,q$, in such a way that they turned out decreasing in the dimension. The theorem, however, does not require such an estimate, and in fact the same estimate is not proved for $w\stp p$. Since one of the estimates for the Reynolds stress requires $L^2_tL^{\xfr{2r}{2-r}}_x$, I will prove $L^p_tL^\8_x$ and $L^2_tL^{\xfr{2r}{2-r}}_x$.
\xend{oss}
\clearpage
\begin{qeddim}
Since the first steps are the same for any exponent pair, we will do them in general, and then deduce the particular cases $L^p_tL^\8_x$ and $L^2_tL^{\xfr{2r}{2-r}}_x$.
\centered*{$L^a_tL^b_x$ estimate}
From the definition, we have:
\[\norm{w\stp c(t)}_{L^b(\T^d)}\leq\sg^{-1}\norm{\sum_k\grad a_k(t)\mathbin:\Wg_k(\sg\per)}_{L^b}\lsim\sg^{-1}\per\!\sum_k\norm{\grad a_k(t)}_{L^\8_x}\norm{\Wg_k(\sg\per)}_{L^b_x}.\]
Now, thanks to \kcref{thm:lemma:Estak}, we take $L^2_t$ to obtain:
\[\norm{w\stp c}_{L^a_tL^b_x}\lsim\sg^{-1}\mg^{\xfr{d-1}{2}-1-\xfr{d-1}{b}}\per\e\sum_k\norm{\grad a_k}_{L^a_tL^\8_x}\leq C_u\sg^{-1}\mg^{\xfr{d-1}{2}-1-\xfr{d-1}{b}}\kg^{\xfr12-\xfr1a}.\]
The desired estimates thus follow from \eqref{eq:Rel3} and \eqref{eq:Rel3b}, the latter of which is obviously satisfied for $\hg$ sufficiently small.
\centered*{$L^s_tW^{\gg,q}_x$ estimate}
This part is very similar to the estimation of $w\stp p$. We first take $W^{\gg,q}_x$ to obtain that:
\begin{align*}
\norm{w\stp c(t)}_{W^{\gg,q}(\T^d)}\leq{}&\sg^{-1}\norm{\sum_k\grad a_k(t)\mathbin:\Wg_k(\sg\per)}_{W^{\gg,q}(\T^d)} \\
{}\lsim{}&\sg^{-1}\per\sum_k\norm{a_k(t)}_{\s{C}^{1+\gg}(\T^d)}\norm{\Wg_k(\sg\per)}_{W^{\gg,q}(\T^d)}.
\end{align*}
Taking $L^s$ norm in time and using \kcref{thm:lemma:Estak} and \kcref{thm:teor:PropMikFlows} we have:
\[\norm{w\stp c}_{L^s_tW^{\gg,q}_x}\lsim\sg^{-1}\per\e\sum_k\norm{a_k(t)}_{L^s_t\s{C}^2_x}\norm{\Wg_k(\sg\per)}_{W^{\gg,q}(\T^d)}\lsim\sg^{\gg-1}\kg^{\xfr12-\xfr1s}\mg^{\gg-1+\xfr{d-1}{2}-\xfr{d-1}{q}}.\]
Once we require \eqref{eq:Rel4}, this implies the desired bound.
\end{qeddim}%
Finally, we estimate the temporal corrector $w\stp t$. From its definition, one can see that the spatial frequency of $w\stp t$ is independent from the parameters $\sg,\tg,\mg$. As a result, this term poses no constraints to the choice of temporal and spatial oscillation/concentration at all and is small for basically any choice of parameters (as long as temporal oscillation $\ng$ is present). This is one of the main technical differences from \cite{BV} and \cite{BCV}, where the leading order effect is temporal oscillation.
\xbegin{propo}[Estimates on the temporal corrector][thm:propo:Estwt]
The temporal corrector $w\stp t$ satisfies:
\[\norm{w\stp t}_{L^\8_tW^{1,\8}_x}\leq C_u\ng^{-1}.\]
In particular, for sufficiently large $\ng$:
\[\norm{w\stp t}_{L^2_tL^2_x}\leq\norm{\lbar R}^{\xfr12}_{L^1_tL^1_x}\qquad\norm{w\stp t}_{L^p_tL^\8_x}+\norm{w\stp t}_{L^s_tW^{\gg,q}_x}\leq\xfr\dg4.\]
\xend{propo}
\begin{qeddim}
It follows directly from the definition of $w\stp t$ that:
\[\norm{w\stp t}_{L^\8_tW^{\gg,\8}_x}\lsim\ng^{-1}\norm{h}_{L^\8_t}\norm{\lbar R}_{L^\8_tW^{\gg+1,\8}_x}\leq C_u\ng^{-1}. \placeqed\]
\end{qeddim}
\clearpage

\ssect{Estimates on the new Reynolds stress}
The last step of the proof is to estimate $R_1$. We proceed with the decomposition in \kcref{thm:lemma:NewReyn}. More specifically, we will prove that, for all sufficiently large $\lg$, each part of the stress $R_1$ is less than $\xfr\dg4$.
\centered*{\tbold{Linear error}}
\xbegin{lemma}[Estimate on the linear error][thm:lemma:EstRlin]
Assume that:
\begin{align*}
2\EXP-1\leq{}&\gg \stag{RelA} \\
r\leq{}&q, \stag{RelB} \\
\ng\kg^{\xfr12}\sg^{-1}\mg^{-1+\xfr{d-1}{2}-\xfr{d-1}{r}}\lsim{}&\lg^{-\hg} \stag{Rel5} \\
\ng^{-1}\lsim\lg^{-\hg}. \stag{RelNu-1}
\end{align*}
Then, for sufficiently large $\lg$,
\[\norm{R_{lin}}_{L^1_tL^r_x}\leq\xfr\dg4.\]
\xend{lemma}
\begin{qeddim}
We split the linear error into three parts:
\[\norm{R_{lin}}_{L^1_tL^r_x}\leq\underbrace{\norm{\s{R}(\xfrl w)}_{L^1_tL^r_x}}_{\eqcolon L_1}+\underbrace{\norm{\s{R}(\pd_t(w\stp p+w\stp c)\!)}_{L^1_tL^r_x}}_{\eqcolon L_2}+\underbrace{\norm{\s{R}(\opn{div}(w\otimes\lbar u+\lbar u\otimes w)\!)}_{L^1_tL^r_x}}_{\eqcolon L_3}.\]
\centered*{Estimate of $L_1$}
By \eqref{eq:AntiDiv1} or boundedness of Riesz transform, we have:
\[L_1\lsim\norm{w}_{L^1_tW^{2\EXP-1,r}_x}.\]
Note that we have estimated $w$ in $L^s_tW^{\gg,q}_x$. Therefore, given \eqref{eq:RelA} and \eqref{eq:RelB}, by the above Propositions we can conclude that:
\[L_1\leq C_u\lg^{-\hg}.\]
\centered*{Estimate of $L_2$}
Since $w\stp p+w\stp c=\sg^{-1}\opn{div}\sum_ka_k(t,x)\Wg_k(\sg x)$, which we saw above, we have:
\[\pd_t(w\stp p+w\stp c)=\sg^{-1}\per\!\sum_k\opn{div}(\pd_ta_k\Wg(\sg\per)\!),\]
and hence:
\[L_2\leq\norm{\s{R}\pd_t(w\stp p+w\stp c)}_{L^1_tL^r_x}\lsim\sg^{-1}\per\!\sum_k\norm{\s{R}\opn{div}(\pd_ta_k\Wg_k(\sg\per)\!)}_{L^1_tL^r_x}.\]
Since $\s{R}\opn{div}$ is a Calderón-Zygmund operator on $\T^d$, we have:
\[L_2\lsim\sg^{-1}\per\!\sum_k\norm{\pd_ta_k}_{L^1_tL^\8_x}\norm{\Wg_k}_{L^r}.\]
Appealing to \kcref{thm:lemma:Estak} and estimates of $\Wg_k$ listed in \kcref{thm:teor:PropMikFlows}, we have:
\[L_2\leq C_u\ng\kg^{\xfr12}\sg^{-1}\mg^{-1+\xfr{d-1}{2}-\xfr{d-1}{r}}.\]
The desired bound then follows from \eqref{eq:Rel5}.
\centered*{Estimate of $L_3$}
For the last term, we simply use the $L^r$ boundedness of $\s{R}$, a crude bound, and the estimates obtained above, in conjunction to \eqref{eq:RelNu-1}, to obtain:
\[L_3\lsim\norm{w\otimes\lbar u}_{L^1_tL^r_x}\lsim\norm{w}_{L^p_tL^\8_x}\norm{\lbar u}_{L^\8_tL^\8_x}\lsim_u\lg^{-\hg}.\]
From these three estimates we can conclude that, for all sufficiently large $\lg$, there holds:
\[\norm{R_{lin}}_{L^1_tL^r_x}\lsim\lg^{-\hg}\leq\xfr\dg4. \placeqed\]
\end{qeddim}
\centered*{\tbold{Correction error}}
\xbegin{lemma}[Estimates on the correction error][thm:lemma:EstRcorr]
For sufficiently large $\lg$:
\[\norm{R_{cor}}_{L^1_tL^r_x}\leq\xfr\dg4.\]
\xend{lemma}
\begin{qeddim}
By boundedness of $\s{R}\opn{div}$ in $L^r$ and Hölder's inequality:
\begin{align*}
\norm{R_{cor}}_{L^1_tL^r_x}\lsim{}&\norm{(w\stp c+w\stp t)\otimes w}_{L^1_tL^r_x}+\norm{w\stp p\otimes(w\stp c+w\stp t)}_{L^1_tL^r_x}\lsim{} \\
{}\lsim{}&\pa{\norm{w\stp c}_{L^2_tL^{\xfr{2r}{2-r}}_x}+\norm{w\stp t}_{L^2_tL^\8_x}}\norm{w}_{L^2_tL^2_x} \\
&{}+\norm{w\stp p}_{L^2_tL^2_x}\pa{\norm{w\stp c}_{L^2_tL^{\xfr{2r}{2-r}}_x}+\norm{w\stp t}_{L^2_tL^\8_x}}.
\end{align*}
Note now that, for $r\in[1,2)$, $\xfr{2r}{2-r}\geq2$, and thus $L^{\xfr{2r}{2-r}}\hra L^2$. Using the propositions above, we can conclude, once again requiring $\ng^{-1}\lsim\lg^{-\hg}$ for the temporal corrector:
\begin{align*}
\norm{w}_{L^2_tL^2_x}\lsim{}&\norm{w\stp p}_{L^2_tL^2_x}+\norm{w\stp c}_{L^2_tL^2_x}+\norm{w\stp t}_{L^2_tL^2_x}\lsim\norm{\lbar R}_{L^1_tL^1_x}^{\xfr12}, \\[0.2em]
\norm{w\stp c}_{L^2_tL^{\xfr{2r}{2-r}}_x}+\norm{w\stp t}_{L^2_tL^\8_x}\leq{}&C_u\lg^{-\hg}.
\end{align*}
This completes the proof. \placeqed*
\end{qeddim}
\centered*{Oscillation error}
The following is our counterpart to \cite[Lemma 5.8]{CL}.
\xbegin{lemma}[Estimates on the oscillation error][thm:lemma:EstRosc]
Assume that:
\begin{align*}
\sg^{-1}\mg^{d-1-\xfr{d-1}{r}}\lsim{}&\lg^{-\hg} \stag{RelOscX} \\
\sg^{-1}\lsim{}&\lg^{-\hg} \stag{RelOscT} \\
\mg^{d-1-\xfr dr}\lsim{}&\lg^{-\hg}. \stag{RelFar}
\end{align*}
Then, for sufficiently large $\lg$,
\[\norm{R_{osc}}_{L^1_tL^r_x}\leq\xfr\dg4.\]
\xend{lemma}
\begin{qeddim}
We will use the decomposition from \kcref{thm:lemma:OscReynPress}:
\[R_{osc}=R_{osc,x}+R_{osc,t}+R_{far}.\]
\centered*{Estimate of $R_{osc,x}$}
As is shown in \cite{CL}, the following estimate holds:
\[\norm{R_{osc,x}}_{L^1_tL^r_x}\leq C_u\sg^{-1}\mg^{d-1-\xfr{d-1}{r}}.\]
\centered*{Estimate of $R_{osc,t}$}
Using the bound on $g_\kg$, we infer:
\[\norm{R_{osc,t}}_{L^1_tL^r_x}=\norm{\sg^{-1}h_\kg(\sg t)\opn{div}\pd_t\lbar R}_{L^1_tL^r_x}\lsim\sg^{-1}\norm{h_\kg(\sg\per)}_{L^1}C_u\leq C_u\sg^{-1}.\]
\centered*{Estimate of $R_{far}$}
We can use \kcref{thm:teor:PropMikFlows} and \kcref{thm:lemma:Estak} to obtain:
\begin{align*}
\norm{R_{far}}_{L^1_tL^r_x}={}&\norm{\sum_{k\neq k'}a_ka_{k'}W_k(\sg\per)\otimes W_{k'}(\sg\per)}_{L^1_tL^r_x}\lsim{} \\
{}\lsim{}&\sum_{k\neq k'}\norm{a_k}_{L^2_tL^\8_x}\norm{a_{k'}}_{L^2_tL^\8_x}\norm{W_k\otimes W_{k'}}_{L^r}\leq{} \\
{}\leq{}&C_u\mg^{d-1-\xfr dr}.
\end{align*}
The desired estimate thus follows from \eqref{eq:RelOscX}-\eqref{eq:RelFar}, completing the proof. \placeqed*
\end{qeddim}

\sect{Choice of parameters} \label{Calcs}
\ssect{General system}
Let us assemble a system of relations required for the above to work:
\begin{align*}
\mg,\ng,\sg\in{}&\N \\
\gg\leq{}&d \tag{\eqref{eq:QG}} \\
q\leq{}&\xfr{2d}{(2\gg-d)^+} \tag{\eqref{eq:QG}} \\
\sg\lsim{}&\ng^2 \tag{\eqref{eq:Rel0}} \\
\kg^{\xfr12-\xfr1p}\mg^{\xfr{d-1}{2}}\lsim{}&\lg^{-\hg} \tag{\eqref{eq:Rel1}} \\
\sg^\gg\kg^{\xfr12-\xfr1s}\mg^{\gg+\xfr{d-1}{2}-\xfr{d-1}{q}}\lsim{}&\lg^{-\hg} \tag{\eqref{eq:Rel2}} \\
\sg^{-1}\mg^{\xfr{d-1}{2}-1}\kg^{\xfr12-\xfr1p}\lsim{}&\lg^{-\hg} \tag{\eqref{eq:Rel3}} \\
\sg^{-1}\mg^{d-2-\xfr{d-1}{r}}\lsim{}&\lg^{-\hg} \tag{\eqref{eq:Rel3b}} \\
\sg^{\gg-1}\kg^{\xfr12-\xfr1s}\mg^{\gg-1+\xfr{d-1}{2}-\xfr{d-1}{q}}\lsim{}&\lg^{-\hg} \tag{\eqref{eq:Rel4}} \\
2\EXP-1\leq{}&\gg \tag{\eqref{eq:RelA}} \\
r\leq{}&q \tag{\eqref{eq:RelB}} \\
\ng\kg^{\xfr12}\sg^{-1}\mg^{-1+\xfr{d-1}{2}-\xfr{d-1}{r}}\lsim{}&\lg^{-\hg} \tag{\eqref{eq:Rel5}} \\
\ng^{-1}\lsim{}&\lg^{-\hg} \tag{\eqref{eq:RelNu-1}} \\
\sg^{-1}\mg^{d-1-\xfr{d-1}{r}}\lsim{}&\lg^{-\hg} \tag{\eqref{eq:RelOscX}} \\
\sg^{-1}\lsim{}&\lg^{-\hg} \tag{\eqref{eq:RelOscT}} \\
\mg^{d-1-\xfr dr}\lsim{}&\lg^{-\hg}. \tag{\eqref{eq:RelFar}}
\end{align*}
We now eliminate a couple of redundancies.
\bi
\item Firstly, if we multiply the LHS of \eqref{eq:Rel4} by $\sg\mg>1$, we  get \eqref{eq:Rel2}. \eqref{eq:Rel4} is therefore a consequence of \eqref{eq:Rel2}, and may be neglected.
\item Similarly, if we multiply the LHS of \eqref{eq:Rel3b} by  $\mg>1$, we obtain \eqref{eq:RelOscX}, which is therefore strictly stronger. We can neglect \eqref{eq:Rel3b}.
\item Continuing, if we multiply the LHS of \eqref{eq:Rel3} by $\sg\mg>1$, we obtain \eqref{eq:Rel1}. Thus, we neglect \eqref{eq:Rel3};
\item Finally, if we multiply the LHS of \eqref{eq:RelOscT} by $\mg^{d-1-\xfr{d-1}{r}}>1$ (recall that $r>1$), we obtain the LHS of \eqref{eq:RelOscX}, meaning \eqref{eq:RelOscT} can be neglected.
\ei
We now choose:
\[(\mg,\ng,\sg,\kg)\coloneq(\lg^\ag,\lg^\bg,\lg^\dg,\lg^\zg).\]
Taking base-$\lg$ logs, the system above becomes:
\begin{align*}
\mg,\ng,\sg\in{}&\N \\
\gg\leq{}&d \tag{\eqref{eq:QG}} \\
q\leq{}&\xfr{2d}{(2\gg-d)^+} \tag{\eqref{eq:QG}} \\
\dg<{}&2\bg \tag{\eqref{eq:Rel0}} \\
\zg\pa{\xfr12-\xfr1p}+\ag\xfr{d-1}{2}<{}&-\hg \tag{\eqref{eq:Rel1}} \\
\dg\gg+\zg\pa{\xfr12-\xfr1s}+\ag\pa{\gg+\xfr{d-1}{2}-\xfr{d-1}{q}}<{}&-\hg \tag{\eqref{eq:Rel2}} \\
2\EXP-1\leq{}&\gg \tag{\eqref{eq:RelA}} \\
r\leq{}&q \tag{\eqref{eq:RelB}} \\
\bg+\xfr\zg2-\dg+\ag\pa{-1+\xfr{d-1}{2}-\xfr{d-1}{r}}<{}&-\hg \tag{\eqref{eq:Rel5}} \\
-\bg<{}&-\hg \tag{\eqref{eq:RelNu-1}} \\
-\dg+\ag\pa{d-1-\xfr{d-1}{r}}<{}&-\hg \tag{\eqref{eq:RelOscX}}  \\
\ag\pa{d-1-\xfr dr}<{}&-\hg. \tag{\eqref{eq:RelFar}} 
\end{align*}
Since $\hg$ is supposed to be small, and $r$ is supposed to be close to 1, we can just substitute $r=1,\hg=0$: all other relations are strict inequalities, which will leave room for an $\hg\ll1$ and a $r\sim1$ to be found so that the original system is satisfied. This means that:
\bi
\item \eqref{eq:RelNu-1} become obvious since $\ag,\bg,\dg,\zg>0$;
\item $r\leq q$ reduces to $q>1$;
\item \eqref{eq:RelOscX} reduces to $-\dg<0$ when we replace $r=1,\hg=0$, so we can remove it;
\item Analogously, \eqref{eq:RelFar} reduces to $-\ag<0$, which is obvious when $\ag>0$.
\ei
We then substitute $(\bg',\dg',\zg')\coloneq\ag^{-1}(\bg,\dg,\zg)$ and eliminate the parameter $\ag>0$. We also note that \eqref{eq:Rel1} immediately prevents $p$ from going above 2, so that we can divide both sides by $\xfr{2-p}{2p}>0$. With simple algebra on top of all this, we can rewrite the system as the system found in \kcref{thm:teor:Relaz}, thus completing the theorem's proof.
\begin{align*}
\mg,\ng,\sg\in{}&\N \\
\bg',\dg',\zg'>{}&0 \\
p,s\geq{}&1 \\
q>{}&1 \\
\dg'<{}&2\bg' \stag{A} \\
\xfr{2p}{2-p}\xfr{d-1}{2}<{}&\zg' \stag{B} \\
sK_{\dg',\gg,\zg',d,q}\coloneq s\pa{\dg'\gg+\xfr{\zg'}{2}+\gg+\xfr{d-1}{2}-\xfr{d-1}{q}}<{}&\zg' \stag{C} \\
2\bg'+\zg'<{}&2\dg'+d+1 \stag{D} \\
2\EXP\leq{}&\gg+1 \stag{E} \\
\gg\leq{}&d \stag{F1} \\
q\leq{}&\xfr{2d}{(2\gg-d)^+}. \stag{F2} \\
\end{align*}
The quantity $K_{\dg',\gg,\zg',d,q}$ in \eqref{eq:C} is nonpositive if
\[q\pa{\dg'\gg+\xfr{\zg'}{2}+\gg+\xfr{d-1}{2}}\leq d-1\iff q\leq q_{min}(d,\dg',\gg,\zg')\coloneq\xfr{2d-2}{2\dg'\gg+\zg'+d-1+2\gg}.\]
In order for this to be possible, we need $q_{min}>1$. However, we see that
\[q_{min}<\xfr{2d-2}{\xfr{2}{2-p}(d-1)+2\gg}<2-p\leq1,\]
due to \eqref{eq:B} and the facts that $\dg'>0,p\geq1$. Thus, $K_{\dg',\gg,\zg',d,q}>0$, and we can divide by it, rewriting \eqref{eq:C} as:
\[s<s_{max}(\dg',\zg'd,q,\gg)\coloneq\xfr{2\zg'}{2\dg'\gg+\zg'+d-1+2\gg-\xfr2q(d-1)}.\]
To have $s_{max}>1$, we require that:
\[\zg'>2\dg'\gg+d-1+2\gg-\xfr2q(d-1)\iff qH_{\dg',\gg,d,\zg'}\coloneq q\pa{2\dg'\gg+d-1+2\gg-\zg'}<2d-2.\]
We remark that
\[2\gg(1+\dg')-\xfr{2p-2}{2-p}(d-1)>H_{\dg',\gg,d,\zg'}>d-1+2\gg-2(\dg'-\bg')-d-1=2(\gg-\dg'+\bg'-1),\]
which unfortunately does not let us determine the sign of $H_{\dg',\gg,d,\zg'}$. If $H_{\dg',\gg,d,\zg'}>0$, we can divide by it, and thus rewrite $s_{max}>1$ as:
\[q<\xfr{2d-2}{2\dg'\gg+d-1+2\gg-\zg'}.\]
Otherwise $s_{max}>1$ for any choice of $q$, so we can write:
\[q<q_{max}(\dg',\zg',d,\gg)\coloneq\xfr{2d-2}{(2\dg'\gg+d-1+2\gg-\zg')^+},\]
where ``<1/0'' is taken to mean ``arbitrary''. To have $q_{max}>1$ as well, we will need either $2\dg'\gg+d+2\gg<\zg'+1$, or:
\[d-1>2\gg(\dg'+1)-\zg'\iff\gg<\xfr{d-1+\zg'}{2(\dg'+1)}\iff\zg'>2\gg(\dg'+1)+1-d.\]
Compatibility of $q_{max}>1$ with \eqref{eq:D} means:
\[2\dg'+d+1-2\bg'>2\gg(\dg'+1)+1-d\iff2d>2\bg'+2\dg'(\gg-1)+2\gg. \stag{comp1} \]
Compatibility of $q_{max}>1$ with $H_{\dg',\gg,d,\zg'}>0$ requires:
\[2\gg(\dg'+1)+1-d<2\dg'\gg+d-1+2\gg\iff2d-2>0,\]
which is of course true.
\nipar If instead $H_{\dg',\gg,d,\zg'}<0$, we will need this relation to be compatible with \eqref{eq:D}, which translates to
\[2\dg'\gg+d-1+2\gg<2(\dg'-\bg')+d+1\iff\gg(\dg'+1)<\dg'-\bg'+1. \stag{comp2}\]
With all of this, we split the system of relations into two cases: $K_{\dg',\gg,\zg',d,q},H_{\dg',\gg,d,\zg'}>0$ and $K_{\dg',\gg,\zg',d,q}>0,H_{\dg',\gg,d,\zg'}<0$. The relation system for the first case reads:
\[\left\{\begin{array}{lr}
\mg,\ng,\sg\in\N \\
\bg',\dg',\zg'>0 \\
p\in[1,2) \\
s\geq1 \\
q>1 \\
\dg'<2\bg' & \text{\eqref{eq:A}} \\
\zg'>\xfr{p}{2-p}(d-1) & \text{\eqref{eq:B}} \\
s<s_{max}(\dg',\zg'd,q,\gg)=\xfr{2\zg'}{2\dg'\gg+\zg'+d-1+2\gg-\xfr2q(d-1)} & \text{\eqref{eq:C}} \\
q<q_{max}(\dg',\zg',d,\gg)=\xfr{2d-2}{2\dg'\gg+d-1+2\gg-\zg'} & s_{max}>1 \\
\zg'<2\dg'\gg+d-1+2\gg & H_{\dg',\gg,d,\zg'}>0 \\
2\gg(\dg'+1)+1-d<\zg' & q_{max}>1 \\
\bg'+\dg'(\gg-1)+\gg<d & \text{\eqref{eq:comp1}} \\
\zg'<2(\dg'-\bg')+d+1 & \text{\eqref{eq:D}} \\
2\EXP\leq\gg+1 \\
\gg\leq d \\
q\leq\xfr{2d}{(2\gg-d)^+}
\end{array}\right. \stag{K+H+} \]
The relation system for the second one reads:
\[\left\{\begin{array}{lr}
\mg,\ng,\sg\in\N \\
\bg',\dg',\zg'>0 \\
p\in[1,2) \\
s\geq1 \\
q>1 \\
\dg'<2\bg' & \text{\eqref{eq:A}} \\
\zg'>\xfr{p}{2-p}(d-1) & \text{\eqref{eq:B}} \\
s<s_{max}(\dg',\zg'd,q,\gg)=\xfr{2\zg'}{2\dg'\gg+\zg'+d-1+2\gg-\xfr2q(d-1)} & \text{\eqref{eq:C}} \\
\zg'\geq2\dg'\gg+d-1+2\gg & H_{\dg',\gg,d,\zg'}\leq0 \\
\gg(\dg'+1)<\dg'-\bg'+1 & \text{\eqref{eq:comp2}} \\
\zg'<2(\dg'-\bg')+d+1 & \text{\eqref{eq:D}} \\
2\EXP\leq\gg+1 \\
\gg\leq d \\
q\leq\xfr{2d}{(2\gg-d)^+}
\end{array}\right. \stag{K+H-} \]
Let us start from the second system, \eqref{eq:K+H-}. $q$ is arbitrary but larger than 1, since it has no upper bounds. As for $s$, we see that $s_{max}$ is $\dg'$-decreasing, so we would like to set $\dg'=0$. With this choice, \eqref{eq:comp2} becomes $\gg<1-\bg'$, which means $\gg\leq d$ and $q\leq\xfr{2d}{(2\gg-d)^+}$ are automatic. We thus rewrite the system:
\[\left\{\begin{array}{lr}
\mg,\ng,\sg\in\N \\
\bg',\zg'>0 \\
p\in[1,2) \\
s\geq1 \\
q>1 \\
\zg'>\xfr{p}{2-p}(d-1) & \text{\eqref{eq:B}} \\
s<s_{max}'(\zg'd,q,\gg)\coloneq\xfr{2\zg'}{\zg'+d-1+2\gg-\xfr2q(d-1)} & \text{\eqref{eq:C}} \\
\zg'\geq2\gg+d-1 & H_{\dg',\gg,d,\zg'}\leq0 \\
\gg<1-\bg' & \text{\eqref{eq:comp2}} \\
\zg'<d+1-2\bg' & \text{\eqref{eq:D}} \\
2\EXP\leq\gg+1 \\
\end{array}\right.\]
$\bg'=0$ is the optimal choice for regularity range. We now test the compatibility of \eqref{eq:D} and \eqref{eq:B}, which with $\bg'=0$ reduces to
\[d+1>\xfr{p}{2-p}(d-1)\iff(2-p)(d+1)>p(d-1)\iff2pd<2d+2\iff p<1+\xfr1d.\]
Let us now compare \eqref{eq:B} with $H<0$. The former will be the stricter bound if and only if
\begin{align*}
\xfr{p}{2-p}(d-1)>d-1+2\gg\iff{}&\gg<\xfr{p-1}{2-p}(d-1) \\
{}\iff{}&2\gg+d-1<p(\gg+d-1) \\
{}\iff{}&p>p_0(\gg,d)\coloneq1+\xfr{\gg}{\gg+d-1}>1.
\end{align*}
We now remark that, since $\gg<1$, $p_0$ satisfies $p<1+\xfr1d$ if and only if
\[\xfr{\gg}{\gg+d-1}<\xfr1d\iff\gg-2+d(1-\gg)>0\iff d>d_0(\gg)\coloneq\xfr{1-\gg}{2-\gg},\]
which is automatic since $d_0<\xfr12<2\leq d$. This means that we have two cases: $p<p_0$, where $H<0$ is stricter than \eqref{eq:B}, and $p_0\leq p<1+\xfr1d$, where the opposite is true.
\nipar With all of this, we wonder if $s_{max}'(\zg',d,q,\gg)\coloneq s_{max}(0,\zg',d,q,\gg)$ is $\zg'$-decreasing. Let us compute the derivative:
\[\pd_{\zg'}s_{max}'(\zg',d,q,\gg)=\xfr{2\pa{1-\xfr2q}(d-1)+4\gg}{\sq{\zg'+2\gg+\pa{1-\xfr2q}(d-1)}^2}.\]
If $q>2$, we can guarantee that this is positive, meaning that $s_{max}'$ is $\zg'$-increasing, and thus is bounded by the value at $d+1$. Assume instead that the numerator is negative. Then $s_{max}$ can do one of two things:
\bi
\item Reach infinity before $\zg'$ reaches its lower bound $\xfr{p}{2-p}(d-1)$, in which case the upper bound on $s_{max}$ will be $s_{max}<\8$ for $\zg'$ somewhere in the range;
\item Not do that, in which case $s_{max}$ is capped by its value at the largest of $\zg'_0(p,d)\coloneq\xfr{p}{2-p}(d-1)$ and $\zg'_1(\gg,d)\coloneq2\gg+d-1$.
\ei
In other words, one of the below will hold:
\begin{align*}
s_{max}'\leq{}&s_{max}'(\max\{\zg'_0,\zg'_1\},d,q,\gg) \\
s_{max}'<{}&\8\quad\wedge\quad s_{max}'(\max\{\zg'_0,\zg'_1\},d,q,\gg)<0.
\end{align*}
We now investigate the possibility of option 2.
\begin{align*}
s_{max}'(\zg'_0,d,q,\gg)={}&\xfr{2\xfr{p}{2-p}(d-1)}{\sq{\xfr{p}{2-p}+1-\xfr2q}(d-1)+2\gg} \\
{}={}&\xfr{2p(d-1)}{\sq{p+(2-p)\pa{1-\xfr2q}}(d-1)+(4-2p)\gg} \\
{}={}&\xfr{p(d-1)}{\sq{1-\xfr1q(2-p)}(d-1)+(2-p)\gg}.
\end{align*}
Set
\[f(p,d,\gg,q)\coloneq d-1+(2-p)\pa{\gg-\xfr{d-1}{q}}.\]
$f$ is clearly monotonic in $p$, whether increasing or decreasing. We can thus bound it as follows:
\begin{align*}
f(p,d,\gg,q)\geq{}&\min\br{f(1,d,\gg,q),f\pa{1+\xfr1d,d,\gg,q}} \\
{}={}&\min\br{d-1+\gg-\xfr{d-1}{q},d-1+\pa{1-\xfr1d}\pa{\gg-\xfr{d-1}{q}}} \\
{}={}&\min\br{(d-1)\pa{1-\xfr1q}+\gg,(d-1)\pa{1-\xfr{d-1}{qd}}+\pa{1-\xfr1d}\gg}>0,
\end{align*}
since $d-1<d<qd,q>1$. So $f\leq0$ is incompatible with the system. This implies that, if $s_{max}$ is negative at $\zg'_1$, then $\zg'_1<\zg'_0$, so the value of $s_{max}$ at $\zg'_1$ does not matter to us. Therefore, we have the following bound on $s$:
\begin{align*}
s<{}&s_{maxmax}(d,q,\gg,p) \\
{}\coloneq{}&\max\br{s_{max}\pa{\zg'_0,d,q,\gg}\,\,,\,\,s_{max}(d+1,d,q,\gg),s_{max}(2\gg+d-1,d,q,\gg)} \\
{}={}&\max\br{\xfr{p(d-1)}{f(p,d,\gg,q)}\,\,,\,\,\xfr{d+1}{d+\gg-\xfr{d-1}{q}},\xfr{2\gg+d-1}{2\gg+d-1-\xfr1q(d-q)}}.
\end{align*}
To write out the ranges of the various cases explicitly, let us see when $\pd_{\zg'}s'_{max}>0$:
\[2\pa{1-\xfr2q}(d-1)+4\gg>0\iff q(4\gg+2d-2)>4(d-1)\iff q>\xfr{2(d-1)}{2\gg+d-1}.\]
The above discussion gives us three ranges:
\[\begin{sistema}
2\EXP\leq\gg+1 \\
\gg<1 \\
p<1+\xfr1d \\
q>\xfr{2d-2}{2\gg+d-1} \\
s<\xfr{d+1}{d+\gg-\xfr{d-1}{q}}
\end{sistema}
\qquad
\begin{sistema}
2\EXP\leq\gg+1 \\
\gg<1 \\
p<1+\xfr{\gg}{\gg+d-1} \\
q\leq\xfr{2d-2}{2\gg+d-1} \\
s<\xfr{2\gg+d-1}{2\gg+d-1-\xfr{d-1}{q}}
\end{sistema}
\qquad
\begin{sistema}
2\EXP\leq\gg+1 \\
\gg<1 \\
p\ni\sqpa{1+\xfr{\gg}{\gg+d-1},1+\xfr1d} \\
q\leq\xfr{2d-2}{2\gg+d-1} \\
s<\xfr{p(d-1)}{d-1+(2-p)\pa{\gg-\xfr{d-1}{q}}}
\end{sistema}.\]
Let us review the last system left to study, namely \eqref{eq:K+H+}:
\[\left\{\begin{array}{lr}
\mg,\ng,\sg\in\N \\
\bg',\dg',\zg'>0 \\
p\in[1,2) \\
s\geq1 \\
q>1 \\
\dg'<2\bg' & \text{\eqref{eq:A}} \\
\zg'>\xfr{p}{2-p}(d-1) & \text{\eqref{eq:B}} \\
s<s_{max}(\dg',\zg'd,q,\gg)=\xfr{2\zg'}{2\dg'\gg+\zg'+d-1+2\gg-\xfr2q(d-1)} & \text{\eqref{eq:C}} \\
q<q_{max}(\dg',\zg',d,\gg)=\xfr{2d-2}{2\dg'\gg+d-1+2\gg-\zg'} & s_{max}>1 \\
\zg'<2\dg'\gg+d-1+2\gg & H_{\dg',\gg,d,\zg'}>0 \\
2\gg(\dg'+1)+1-d<\zg' & q_{max}>1 \\
\bg'+\dg'(\gg-1)+\gg<d & \text{\eqref{eq:comp1}} \\
\zg'<2(\dg'-\bg')+d+1 & \text{\eqref{eq:D}} \\
2\EXP\leq\gg+1 \\
\gg\leq d \\
q\leq\xfr{2d}{(2\gg-d)^+}
\end{array}\right. \]
The compatibility of $H>0$ and \eqref{eq:B} requires
\[2\dg'\gg+d-1+2\gg>\xfr{p}{2-p}(d-1)\iff\gg(1+\dg')>\xfr{p-1}{2-p}(d-1).\]
This is compatible with $\gg\leq d$ iff
\[d>\xfr{p-1}{2-p}\xfr{d-1}{1+\dg'}\iff\dg'd>\xfr{p-1}{2-p}(d-1)-d=\xfr{2pd-3d-p+1}{2-p}.\]
Since both $q_{max}$ and $s_{max}$ are $\dg'$-decreasing, we would like to take $\dg'=0$, which combined with this condition would mean:
\[2pd-3d-p+1\leq0\iff d(2p-3)\leq p-1.\]
Here we either assume $p\leq\xfr32$, in which case this is obvious, or this is not true. So these compatibilities reduce to $p\leq\xfr32$. Once we take $\dg'=0$, we will also want to take $\bg'=0$, so as to maximize the ranges of both $\zg'$ and $\gg$. This will lead us to $\zg'<d+1$, which combined with \eqref{eq:B} leads to $(d+1)(2-p)>p(d-1)$, or $2d+2-p>2pd-p$, or $p<1+\xfr1d$, which implies $p<\xfr32$.
\nipar In fact, $p<1+\xfr1d$ implies the compatibility above follows by assuming $\gg\geq1$. Under that assumption, we compare \eqref{eq:D} and $H>0$, and see that $d+1<2\gg+d-1\iff2\gg-2>0$, which holds. This means that the stricter upper bound is always \eqref{eq:D}, and we can thus neglect $H>0$.
\nipar With all of that, we rewrite the system:
\[\left\{\begin{array}{lr}
\mg,\ng,\sg\in\N \\
\zg'>0 \\
p\in\sqpa{1,1+\xfr1d} \\
s\geq1 \\
q>1 \\
\zg'>\xfr{p}{2-p}(d-1) & \text{\eqref{eq:B}} \\
s<s_{max}(0,\zg'd,q,\gg)=\xfr{2\zg'}{\zg'+d-1+2\gg-\xfr2q(d-1)}\eqcolon s_{max}'(\zg',d,q,\gg) & \text{\eqref{eq:C}} \\
q<q_{max}(0,\zg',d,\gg)=\xfr{2d-2}{d-1+2\gg-\zg'}\eqcolon q_{max}'(\zg',d,\gg) & s_{max}>1 \\
2\gg+1-d<\zg' & q_{max}>1 \\
\zg'<d+1 & \text{\eqref{eq:D}} \\
\gg\in[1,d) \\
2\EXP\leq\gg+1 \\
q\leq\xfr{2d}{(2\gg-d)^+}
\end{array}\right. \]
Since $q_{max}'$ is clearly $\zg'$-increasing, the optimal choice for it would be to maximize $\zg'$, i.e. choose $\zg'=d+1$. Let us now investigate the monotonicity of $s_{max}'$ w.r.t. $\zg'$:
\begin{align*}
\pd_{\zg'}s_{max}'(\zg',d,q,\gg)={}&\xfr{2\zg'+2d-2+4\gg-\xfr4q(d-1)-2\zg'}{\pa{\zg'+d-1+2\gg-\xfr2q(d-1)}^2}={} \\
{}={}&\xfr{\xfr{d-1}{2}+\gg-\xfr1q(d-1)}{\pa{\xfr{\zg'+d-1}{2}+\gg-\xfr1q(d-1)}^2}={} \\
{}={}&\xfr{\gg+\pa{\xfr12-\xfr1q}(d-1)}{\pa{\xfr{\zg'+d-1}{2}+\gg-\xfr1q(d-1)}^2}.
\end{align*}
The condition for this to be non-positive is
\[\gg+\pa{\xfr12-\xfr1q}(d-1)\leq0\iff q\leq\xfr{2d-2}{2\gg+d-1}=2-\xfr{4\gg}{2\gg+d-1}\eqcolon q_0(d,\gg)=q_{max}'(0,d,\gg).\]
$q_0<q_{max}$ whenever $\zg'>0$, which is a given. $q_0\leq2$ whenever $\gg\geq0$, also a given. $q_0>1$ iff $d-1>2\gg$, or $\gg<\xfr{d-1}{2}$, which is allowed. So for large $\gg$ the system will not allow $q\leq q_0$, and thus $s_{max}'$ will be $\zg'$-increasing. For small enough $\gg$, however, if $q$ is large enough, $s_{max}$ is $\zg'$-increasing, otherwise it is $\zg'$-decreasing. For $q=q_0$ we have $s_{max}\equiv2$.
\nipar In any case, the sign of this derivative is independent of $\zg'$, so the monotonicity is guaranteed. We have three bounds on $\zg'$: the two lower bounds \eqref{eq:B} and $q_{max}>1$, and the one upper bound \eqref{eq:D}. Concerning \eqref{eq:B} and $q_{max}>1$, let us investigate the relation between their bounding quantities.
\begin{align*}
2\gg+1-d>\xfr{p}{2-p}(d-1)\iff{}&\gg>\xfr{d-1}{2-p}\iff2-p>\xfr{d-1}{\gg} \\
{}\iff{}&p<2-\xfr{d-1}{\gg}=\xfr{2\gg+1-d}{\gg}.
\end{align*}
This means that, in order for \eqref{eq:B} to follow from $q_{max}>1$, we need this quantity to be greater than 1, which boils down to $\gg>d-1$. This is not at all a guarantee, and is in fact close to saturating $\gg<d$. We can therefore conclude the following:
\begin{align*}
s<{}&s_{maxmax}(d,q,p,\gg) \\
{}\coloneq{}&\max\br{s_{max}'\pa{\xfr{p}{2-p}(d-1),d,q,\gg},s_{max}'\pa{2\gg+1-d,d,q,\gg},s_{max}'(d+1,d,q,\gg)}.
\end{align*}
But naturally the choice of $\zg'$ has an impact on $q_{max}$, so we split this into three systems:
\bi
\item $\zg'=d+1,q>q_0$;
\item $\zg'=\xfr{p}{2-p}(d-1),q\leq q_0,\gg\leq d-1$;
\item $\zg'=2\gg+1-d,q\leq q_0,\gg>d-1$.
\ei
However, $q\leq q_0$ requires $q_0>1$ which implies $\gg<\xfr{d-1}{2}$, so we immediately discard the last system. The first system will therefore read:
\[\left\{\begin{array}{lr}
\mg,\ng,\sg\in\N \\
p\in\sqpa{1,1+\xfr1d} \\
s\geq1 \\
q>q_0 \\
s<s_{max}'(d+1,d,q,\gg)=\xfr{d+1}{d+\gg-\xfr1q(d-1)}\eqcolon s_{top,1}(d,q,\gg) \\
q<q_{max}'(d+1,d,\gg)=\xfr{d-1}{\gg-1}\eqcolon q_{top,1}(d,\gg) \\
\gg\in[1,d) \\
2\EXP\leq\gg+1 \\
q\leq\xfr{2d}{(2\gg-d)^+}
\end{array}\right. \stag{Topzg'}\]
For this to yield solutions, we must have $q_{top,2}>q_0$, i.e.:
\[\xfr{d-1}{\gg-1}>\xfr{2d-2}{2\gg+d-1}\iff2\gg+d-1>2\gg-2\iff d>-1,\]
which is true since $d\geq2$.
\nipar We now compare the two upper bounds on $q$. The stricter one is $q\leq\xfr{2d}{(2\gg-d)^+}$ iff $\gg>\xfr d2$ and
\[\xfr{2d}{2\gg-d}\leq\xfr{d-1}{\gg-1}\iff2d(\gg-1)\leq(d-1)(2\gg-d)\iff3d-d^2-2\gg\geq0.\]
Since $\gg>\xfr d2$, we see $3d-d^2-2\gg<2d-d^2=d(2-d)<0$, so that the above never happens. In other words, $q\leq\xfr{2d}{(2\gg-d)^+}$ can be neglected. We thus obtain the following range:
\[\left\{\begin{array}{lr}
\mg,\ng,\sg\in\N \\
p\in\sqpa{1,1+\xfr1d} \\[.4em]
s\in\sqpa{1,\xfr{d+1}{d+\gg-\xfr1q(d-1)}} \\[.8em]
q\in\pa{\xfr{2d-2}{2\gg+d-1},\xfr{d-1}{\gg-1}} \\
\gg\in[1,d) \\
2\EXP\leq\gg+1 \\
\end{array}\right..\]
For the second case, we recall that $q\leq q_0<2$. This automatically guarantees $q\leq\xfr{2d}{(2\gg-d)^+}$, since this is clearly $\gg$-non-increasing and for $\gg=d$ it equals 2. Moreover, we already pointed out that $q_0<q_{max}$, so $q<q_{max}$ is also automatic. One can actually see that $q_{max}>2$ whenever $\gg<\xfr{d-1}{2}$. The system therefore reads:
\[\left\{\begin{array}{lr}
\mg,\ng,\sg\in\N \\
p\in\sqpa{1,1+\xfr1d} \\
s\geq1 \\
1<q\leq q_0(d,\gg)=\xfr{2d-2}{2\gg+d-1}=2-\xfr{4\gg}{2\gg+d-1} \\
s<s_{max}'\pa{\xfr{p}{2-p}(d-1),d,q,\gg}=\xfr{2p(d-1)}{p(d-1)+(2-p)(d-1+2\gg-\xfr2q(d-1)\!)}\eqcolon s_{top,2}(p,d,q,\gg) \\
\gg<\xfr{d-1}{2} \\
2\EXP\leq\gg+1
\end{array}\right. \stag{Lowzg'Lowgg}\]
From the two expressions for $q_0$, we can see that it is $d$-increasing and $\gg$-decreasing, as is expected. $s_{top,2}$ is clearly $\gg$-decreasing. Let us see that it is $d$-increasing:
\begin{align*}
\pd_ds_{top,2}(p,d,q,\gg)={}&\xfr{2p\sq{2d-2+2(2-p)\pa{\gg-\xfr{d-1}{q}}}-2p(d-1)\pa{2-\xfr2q(2-p)}}{\sq{2d-2+2(2-p)\pa{\gg-\xfr{d-1}{q})}}^2} \\
{}={}&\xfr{p(2-p)\pa{\gg-\xfr{d-1}{q}}+p(d-1)\pa{\xfr1q(2-p)}}{\sq{d-1+(2-p)\pa{\gg-\xfr{d-1}{q})}}^2} \\
{}={}&\xfr{p\gg(2-p)}{\sq{d-1+(2-p)\pa{\gg-\xfr{d-1}{q})}}^2}>0.
\end{align*}
We now investigate its behaviour w.r.t. $p$:
\begin{align*}
\pd_ps_{top,2}(p,d,q,\gg)={}&\xfr{4(d-1)\sq{d-1+(2-p)\pa{\gg-\xfr{d-1}{q}}}+4p(d-1)\pa{\gg-\xfr{d-1}{q}}}{\sq{2d-2+2(2-p)\pa{\gg-\xfr{d-1}{q}}}^2} \\
{}={}&\xfr{(d-1)\sq{d-1+2\pa{\gg-\xfr{d-1}{q}}}}{\sq{d-1+(2-p)\pa{\gg-\xfr{d-1}{q}}}^2}\leq0,
\end{align*}
because the square bracket in the numerator is non-positive iff $q\leq q_0$, which is the case we are studying. We therefore have the following ranges:
\[\left\{\begin{array}{lr}
\mg,\ng,\sg\in\N \\
p\in\sqpa{1,1+\xfr1d} \\[.4em]
s\in\sqpa{1,\xfr{p(d-1)}{d-1+(2-p)(\gg-\xfr1q(d-1)\!)}}\sbse\sqpa{1,\xfr{d-1}{(d-1)(1-\xfr1q)+\gg(2-p)}} \\[.8em]
q\in\pasq{1,\xfr{2d-2}{2\gg+d-1}}\sbse(1,2] \\[.4em]
\gg<\xfr{d-1}{2} \\
2\EXP\leq\gg+1
\end{array}\right..\]

\appendix
\sect{Improved Hölder inequality}
In \cite[Appendix B]{CL}, an improvement of the Hölder inequality is stated. Below, we state and prove a generalization of that result.
\xbegin{lemma}[Improved Hölder inequality][thm:lemma:ImprIHI]
Let $\lg\in\N$ and $f,g:\T^d\to\R$ be $\s{C}^s$ function. Let $g_\lg(x)\coloneq g(\lg x)$. Then for every $p\in[1,2],s\in(0,1)$:
\[\abs{\norm{fg_\lg}_{L^p}-\norm{f}_{L^p}\norm{g}_{L^p}}\leq\xfr{C_p}{\lg^{\xfr{s}{p'}}}\norm{f}_{\s{C}^s}\norm{g}_{L^p},\]
where all the norms are taken on $\T^d$ and $p^{-1}+p'^{-1}=1$. In particular:
\[\norm{fg_\lg}_{L^p}\leq\norm{f}_{L^p}\norm{g}_{L^p}+\xfr{C_p}{\lg^{\xfr{s}{p'}}}\norm{f}_{\s{C}^s}\norm{g}_{L^p}.\]
\xend{lemma}
\begin{qeddim}
Let us divide $\T^d$ into $\lg^d$ small cubes $\{Q_j\}_j$ of edge $\xfr1\lg$. On each $Q_j$ we have:
\begin{align*}
\xints{Q_j}{}\2\abs{f(x)}^p\abs{g_\lg(x)}^p\diff x={}&\xints{Q_j}{}\2\pa{\abs{f(x)}^p-\fint_{Q_j}\abs{f(y)}^p\diff y}\abs{g_\lg(x)}^p\diff x \\
&{}+\fint_{Q_j}\abs{f(y)}^p\diff y\xints{Q_j}{}\2\abs{g_\lg(x)}^p\diff x={} \\
{}\pux[\Big]{$g_\lg(x)=g(\lg x)$ and change variables in the integral in term 2}{=}{}&\xints{Q_j}{}\2\pa{\abs{f(x)}^p-\fint_{Q_j}\abs{f(y)}^p\diff y}\abs{g_\lg(x)}^p\diff x \\
&{}+\xfr{1}{\lg^d}\fint{Q_j}\abs{f(y)}^p\diff y\xints{\T^d}{}\2\abs{g(x)}^p\diff x={} \\
{}\pux[\Big]{$|Q_j|=\lg^{-d}$}{=}{}&\xints{Q_j}{}\2\pa{\abs{f(x)}^p-\fint{Q_j}\abs{f(y)}^p\diff y}\abs{g_\lg(x)}^p\diff x \\
&{}+\xints{Q_j}{}\abs{f(y)}^p\diff y\per\e\xints{Q_j}{}\2\abs{g(x)}^p\diff x.
\end{align*}
Summing over $j$ we get:
\[\norm{fg_\lg}_{L^p(\T^d)}^p=\norm{f}_{L^p}^p\norm{g}_{L^p}^p+\sum_j\xints{Q_j}{}\2\pa{\abs{f(x)}^p-\fint_{Q_j}\abs{f(y)}^p\diff y}\abs{g_\lg(x)}^p\diff x.\]
Let us now estimate the second term in the RHS. For $x,y\in Q_j$, we have:
\[\abs[b]{\abs{f(x)}^p-\abs{f(y)}^p}\leq\abs{f(x)}^{p-1}\abs[b]{\abs{f(x)}-\abs{f(y)}}+\abs{f(y}\abs[b]{\abs{f(x)}^{p-1}-\abs{f(y)}^{p-1}}.\]
Since $p-1<1$, raising to the power $p-1$ is a concave function, meaning $||f(x)|^{p-1}-|f(y)|^{p-1}|\lsim|f(x)-f(y)|^{p-1}$. Therefore:
\begin{align*}
\abs[b]{\abs{f(x)}^p-\abs{f(y)}^p}\leq{}&\norm{f}_{\s{C}^0}^{p-1}\norm{f}_{\s{C}^s}\xfr{1}{\lg^s}+K_p\norm{f}_{\s{C}^0}\norm{f}_{\s{C}^s}^{p-1}\xfr{1}{\lg^{s(p-1)}} \\
{}\leq{}&\lg^{-s(p-1)}(K_p+\lg^{-s(2-p)})\norm{f}_{\s{C}^{\xfr12}}^p.
\end{align*}
Therefore:
\begin{align*}
\sum_j\xints{Q_j}{}\2\pa{\abs{f(x)}^p-\fint_{Q_j}\abs{f(y)}^p\diff y}\abs{g_\lg(x)}^p\diff x\leq{}&\xfr{K'_p\norm{f}_{\s{C}^s}^p}{\lg^{s(p-1)}}\per\sum_j\xints{Q_j}{}\2\abs{g_\lg(x)}^p\diff x={} \\
{}={}&\xfr{K'_p\norm{f}_{\s{C}^s}^p}{\lg^{s(p-1)}}\norm{g_\lg}_{L^p}^p={} \\
{}={}&\xfr{K'_p}{\lg^{s(p-1)}}\norm{f}_{\s{C}^s}^p\norm{g}_{L^p}^p.
\end{align*}
This yields:
\[\abs{\norm{fg_\lg}_{L^p}-\norm{f}_{L^p}\norm{g}_{L^p}}^p\lsim_p\abs{\norm{fg_\lg}_{L^p}^p-\norm{f}_{L^p}^p\norm{g}_{L^p}^p}\leq\xfr{K_p'^p}{\lg^{\xfr{sp}{p'}}}\norm{f}_{\s{C}^s}^p\norm{g}_{L^p}^p,\]
where in the first step we used that $p>1$, thus the power $p$ is a convex function. We then take $p$th roots, and obtain the first estimate, where $C_p=K'_p\per H_p$, $H_p$ being the implicit constant in the first step. The second one follows trivially from the first. \placeqed*
\end{qeddim}
\xbegin{oss}[Greater exponents][thm:oss:GreaterExp]
For $p>2$, one gets $p-\lfloor p\rfloor$ instead of $p-1$, and has to insert $\lfloor p\rfloor$ intermediate terms when estimating $||f(x)|^p-|f(y)|^p|$.
\end{oss}

\sect{Antidivergences}
This appendix is devoted to the results on antidivergence operators found in \cite[Appendix B]{CL}.
\nipar For any $f\in\Cinf(\T^d)$, there exists $v\in\Cinf_0(\T^d)$ such that:
\[\Dg v=f-\fint_{\T^d}f.\]
We denote $v$ by $\Dg^{-1}f$. Note that if $f\in\Cinf_t(\T^d)$, then by rescaling we have:
\[\Dg^{-1}(f(\sg\per)\!)=\sg^{-2}v(\sg\per)\qquad\sg\in\N.\]
We recall the following antidivergence operator $\s{R}$.
\xbegin{defi}[Tensor antidivergence][thm:defi:TensAntidiv]
We define the operator $\s{R}:\Cinf(\T^d,\R^d)\to\Cinf(\T^d,\s{S}_0^{d\x d})$ as:
\[(\s{R}v)_{ij}=\s{R}_{ijk}v_k,\]
where:
\[\s{R}_{ijk}=\xfr{d-2}{d-1}\Dg^{-2}\pd_i\pd_j\pd_k-\xfr{1}{d-1}\Dg^{-1}\pd_k\dg_{ij}+\Dg^{-1}\pd_i\dg_{jk}+\Dg^{-1}\pd_j\dg_{ik}.\]
\xend{defi}
It is clear that $\s{R}$ is well defined since $\s{R}_{ijk}$ is symmetric in $i,j$ and taking trace gives:
\begin{align*}
\opn{Tr}\s{R}v={}&-\xfr{d-2}{d-1}\Dg^{-1}\pd_kv_k-\xfr{d}{d-1}\Dg^{-1}\pd_kv_k+\Dg^{-1}\pd_kv_k+\Dg^{-1}\pd_kv_k={} \\
{}={}&\pa{\xfr{2-d}{d-1}-\xfr{d}{d-1}+2}\Dg^{-1}\pd_kv_k=0.
\end{align*}
By a direct computation, one can also show that:
\begin{align*}
\opn{div}(\s{R}v)={}&v-\fint_{\T^d}v\qquad\VA v\in\Cinf(\T^d,\R^d) \\
\opn{div}(\s{R}\Dg v)={}&\grad v+\grad^Tv\qquad\VA v\in\Cinf(\T^d,\R^d):\opn{div}v=0. \tag{B.3}\label{eq:AntiDiv1}
\end{align*}
We can show that $\s{R}$ is bounded in $L^p(\T^d)$ for $1\leq p\leq\8$.
\xbegin{teor}[Lebesgue-space boundedness of the antidivergence][thm:teor:LpBoundR]
Let $1\leq p\leq\8$. For any $f\in\Cinf_0(\T^d)$ there holds:
\[\norm{\s{R}f}_{L^p(\T^d)}\lsim\norm{f}_{\T^d}.\]
\xend{teor}
The proof is left to the paper.
\nipar We can also introduce the bilinear version $\s{B}:\Cinf(\T^d,\R^d)\x\Cinf(\T^d,\s{S}_0^{d\x d}\to\Cinf_0(\T^d,\s{D}_0^{d\x d})$ of $\s{R}$. Let:
\[(\s{B}(v,A)\!)_{ij}\coloneq v_\ell\s{R}_{ijk}A_{\ell k}-\s{R}(\pd_iv_\ell\s{R}_{ijk}A_{\ell k}),\]
or by a slight abuse of notations:
\[\s{B}(v,A)=\s{R}Av-\s{R}(\grad v\s{R}A).\]
\xbegin{teor}[Bilinear antidivergence equation][thm:teor:BilinAntidivEq]
Let $1\leq p\leq\8$. For any $v\in\Cinf(\T^d),A\in\Cinf_0(\T^d,\R^d)$:
\begin{align*}
\opn{div}(\s{B}(v,A)\!)={}&vA-\fint_{\T^d}vA \tag{B.4}\label{eq:Antidiv2} \\
\norm{\s{B}(v,A)}_{L^p(\T^d)}\lsim{}&\norm{v}_{\s{C}^1(\T^d)}\norm{A}_{L^p}(\T^d).
\end{align*}
\xend{teor}
\begin{qeddim}
A direct computation gives:
\begin{align*}
\opn{div}(\s{B}(v,A)\!)={}&\pd_jv_\ell\s{R}_{ijk}A_{\ell k}+v_\ell\pd_j\s{R}_{ijk}A_{\ell k}-\opn{div}\s{R}(\pd_iv_\ell\s{R}_{ijk}A_{\ell k}={} \\
{}={}&v_\ell A_{i\ell}+\fint\pd_iv_\ell\s{R}_{ijk}A_{\ell k},
\end{align*}
where we have used the fact that $A$ has zero mean and $\s{R}$ is symmetric. Integrating by parts, we have:
\[\fint\pd_iv_\ell\s{R}_{ijk}A_{\ell k}=-\fint v_\ell\pd_i\s{R}_{ijk}A_{\ell k}=-\fint v_\ell A_{\ell j},\]
which implies that:
\[\opn{div}(\s{B}(v,A)\!)=vA-\fint vA.\]
Proving \kcref{thm:teor:LpBoundR} involved proving the $\s{R}_{ijk}$ are bounded in $L^p$, so that:
\begin{align*}
\norm{\s{B}(v,A)}_{L^p}\leq{}&\norm{(v_\ell\s{R}_{ijk}A_{\ell k})_{ij}}_{L^p}-\norm{\s{R}(\pd_iv_\ell\s{R}_{ijk}A_{\ell k})}_{L^p}\lsim\norm{v}_{L^p}\norm{A}_{L^p}+\norm{\grad v}_{L^p}\norm{A}_{L^p}\lsim{} \\
{}\lsim{}&\norm{v}_{\s{C}^1}\norm{A}_{L^p}.
\end{align*}
The proof is thus complete. \placeqed*
\end{qeddim}

\end{document}